\input amstex
%

\def\next{AMS-SEKR}\ifx\styname\next \endinput\fi
\catcode`\@=11
\def\styname{AMS-SEKR}
\def\styversion{2.0}
{\W@{}\W@{\styname.STY - Version \styversion}\W@{}}
\hyphenation{acad-e-my acad-e-mies af-ter-thought anom-aly anom-alies
an-ti-deriv-a-tive an-tin-o-my an-tin-o-mies apoth-e-o-ses apoth-e-o-sis
ap-pen-dix ar-che-typ-al as-sign-a-ble as-sist-ant-ship as-ymp-tot-ic
asyn-chro-nous at-trib-uted at-trib-ut-able bank-rupt bank-rupt-cy
bi-dif-fer-en-tial blue-print busier busiest cat-a-stroph-ic
cat-a-stroph-i-cally con-gress cross-hatched data-base de-fin-i-tive
de-riv-a-tive dis-trib-ute dri-ver dri-vers eco-nom-ics econ-o-mist
elit-ist equi-vari-ant ex-quis-ite ex-tra-or-di-nary flow-chart
for-mi-da-ble forth-right friv-o-lous ge-o-des-ic ge-o-det-ic geo-met-ric
griev-ance griev-ous griev-ous-ly hexa-dec-i-mal ho-lo-no-my ho-mo-thetic
ideals idio-syn-crasy in-fin-ite-ly in-fin-i-tes-i-mal ir-rev-o-ca-ble
key-stroke lam-en-ta-ble light-weight mal-a-prop-ism man-u-script
mar-gin-al meta-bol-ic me-tab-o-lism meta-lan-guage me-trop-o-lis
met-ro-pol-i-tan mi-nut-est mol-e-cule mono-chrome mono-pole mo-nop-oly
mono-spline mo-not-o-nous mul-ti-fac-eted mul-ti-plic-able non-euclid-ean
non-iso-mor-phic non-smooth par-a-digm par-a-bol-ic pa-rab-o-loid
pa-ram-e-trize para-mount pen-ta-gon phe-nom-e-non post-script pre-am-ble
pro-ce-dur-al pro-hib-i-tive pro-hib-i-tive-ly pseu-do-dif-fer-en-tial
pseu-do-fi-nite pseu-do-nym qua-drat-ics quad-ra-ture qua-si-smooth
qua-si-sta-tion-ary qua-si-tri-an-gu-lar quin-tes-sence quin-tes-sen-tial
re-arrange-ment rec-tan-gle ret-ri-bu-tion retro-fit retro-fit-ted
right-eous right-eous-ness ro-bot ro-bot-ics sched-ul-ing se-mes-ter
semi-def-i-nite semi-ho-mo-thet-ic set-up se-vere-ly side-step sov-er-eign
spe-cious spher-oid spher-oid-al star-tling star-tling-ly
sta-tis-tics sto-chas-tic straight-est strange-ness strat-a-gem strong-hold
sum-ma-ble symp-to-matic syn-chro-nous topo-graph-i-cal tra-vers-a-ble
tra-ver-sal tra-ver-sals treach-ery turn-around un-at-tached un-err-ing-ly
white-space wide-spread wing-spread wretch-ed wretch-ed-ly Brown-ian
Eng-lish Euler-ian Feb-ru-ary Gauss-ian Grothen-dieck Hamil-ton-ian
Her-mit-ian Jan-u-ary Japan-ese Kor-te-weg Le-gendre Lip-schitz
Lip-schitz-ian Mar-kov-ian Noe-ther-ian No-vem-ber Rie-mann-ian
Schwarz-schild Sep-tem-ber
form per-iods Uni-ver-si-ty cri-ti-sism for-ma-lism}
\Invalid@\nofrills
\Invalid@\usualspace
\newif\ifnofrills@
\def\nofrills@#1#2{\relaxnext@
  \DN@{\ifx\next\nofrills
    \nofrills@true\let#2\relax\DN@\nofrills{\nextii@}%
  \else
    \nofrills@false\def#2{#1}\let\next@\nextii@\fi
\next@}}
\def\usualspace@#1{\ifnofrills@\def\usualspace{#1}\fi}
\def\addto#1#2{\csname \expandafter\eat@\string#1@\endcsname
  \expandafter{\the\csname \expandafter\eat@\string#1@\endcsname#2}}
\newdimen\bigsize@
\def\big@#1#2{{\hbox{$\left#2\vcenter to#1\bigsize@{}%
  \right.\nulldelimiterspace\z@\m@th$}}}
\def\big{\big@\@ne}
\def\Big{\big@{1.5}}
\def\bigg{\big@\tw@}
\def\Bigg{\big@{2.5}}
\def\raggedcenter@{\leftskip\z@ plus.4\hsize \rightskip\leftskip
 \parfillskip\z@ \parindent\z@ \spaceskip.3333em \xspaceskip.5em
 \pretolerance9999\tolerance9999 \exhyphenpenalty\@M
 \hyphenpenalty\@M \let\\\linebreak}
\def\upperspecialchars{\def\ss{SS}\let\i=I\let\j=J\let\ae\AE\let\oe\OE
  \let\o\O\let\aa\AA\let\l\L}
\def\uppercasetext@#1{%
  {\spaceskip1.2\fontdimen2\the\font plus1.2\fontdimen3\the\font
   \upperspecialchars\uctext@#1$\m@th\aftergroup\eat@$}}
\def\uctext@#1$#2${\endash@#1-\endash@$#2$\uctext@}
\def\endash@#1-#2\endash@{\uppercase{#1}\if\notempty{#2}--\endash@#2\endash@\fi}
\def\runaway@#1{\DN@{#1}\ifx\envir@\next@
  \Err@{You seem to have a missing or misspelled \string\end#1 ...}%
  \let\envir@\empty\fi}
\newif\iftemp@
\def\notempty#1{TT\fi\def\test@{#1}\ifx\test@\empty\temp@false
  \else\temp@true\fi \iftemp@}
\font@\tensmc=cmcsc10
\font@\sevenex=cmex7
\font@\sevenit=cmti7
\font@\eightrm=cmr8 
\font@\sixrm=cmr6 
\font@\eighti=cmmi8     \skewchar\eighti='177 
\font@\sixi=cmmi6       \skewchar\sixi='177   
\font@\eightsy=cmsy8    \skewchar\eightsy='60 
\font@\sixsy=cmsy6      \skewchar\sixsy='60   
\font@\eightex=cmex8
\font@\eightbf=cmbx8 
\font@\sixbf=cmbx6   
\font@\eightit=cmti8 
\font@\eightsl=cmsl8 
\font@\eightsmc=cmcsc8
\font@\eighttt=cmtt8 


\loadmsam
\loadmsbm
\loadeufm
\UseAMSsymbols
\newtoks\tenpoint@
\def\tenpoint{\normalbaselineskip12\p@
 \abovedisplayskip12\p@ plus3\p@ minus9\p@
 \belowdisplayskip\abovedisplayskip
 \abovedisplayshortskip\z@ plus3\p@
 \belowdisplayshortskip7\p@ plus3\p@ minus4\p@
 \textonlyfont@\rm\tenrm \textonlyfont@\it\tenit
 \textonlyfont@\sl\tensl \textonlyfont@\bf\tenbf
 \textonlyfont@\smc\tensmc \textonlyfont@\tt\tentt
 \textonlyfont@\bsmc\tenbsmc
 \ifsyntax@ \def\big##1{{\hbox{$\left##1\right.$}}}%
  \let\Big\big \let\bigg\big \let\Bigg\big
 \else
  \textfont\z@=\tenrm  \scriptfont\z@=\sevenrm  \scriptscriptfont\z@=\fiverm
  \textfont\@ne=\teni  \scriptfont\@ne=\seveni  \scriptscriptfont\@ne=\fivei
  \textfont\tw@=\tensy \scriptfont\tw@=\sevensy \scriptscriptfont\tw@=\fivesy
  \textfont\thr@@=\tenex \scriptfont\thr@@=\sevenex
        \scriptscriptfont\thr@@=\sevenex
  \textfont\itfam=\tenit \scriptfont\itfam=\sevenit
        \scriptscriptfont\itfam=\sevenit
  \textfont\bffam=\tenbf \scriptfont\bffam=\sevenbf
        \scriptscriptfont\bffam=\fivebf
  \setbox\strutbox\hbox{\vrule height8.5\p@ depth3.5\p@ width\z@}%
  \setbox\strutbox@\hbox{\lower.5\normallineskiplimit\vbox{%
        \kern-\normallineskiplimit\copy\strutbox}}%
 \setbox\z@\vbox{\hbox{$($}\kern\z@}\bigsize@=1.2\ht\z@
 \fi
 \normalbaselines\rm\ex@.2326ex\jot3\ex@\the\tenpoint@}
\newtoks\eightpoint@
\def\eightpoint{\normalbaselineskip10\p@
 \abovedisplayskip10\p@ plus2.4\p@ minus7.2\p@
 \belowdisplayskip\abovedisplayskip
 \abovedisplayshortskip\z@ plus2.4\p@
 \belowdisplayshortskip5.6\p@ plus2.4\p@ minus3.2\p@
 \textonlyfont@\rm\eightrm \textonlyfont@\it\eightit
 \textonlyfont@\sl\eightsl \textonlyfont@\bf\eightbf
 \textonlyfont@\smc\eightsmc \textonlyfont@\tt\eighttt
 \textonlyfont@\bsmc\eightbsmc
 \ifsyntax@\def\big##1{{\hbox{$\left##1\right.$}}}%
  \let\Big\big \let\bigg\big \let\Bigg\big
 \else
  \textfont\z@=\eightrm \scriptfont\z@=\sixrm \scriptscriptfont\z@=\fiverm
  \textfont\@ne=\eighti \scriptfont\@ne=\sixi \scriptscriptfont\@ne=\fivei
  \textfont\tw@=\eightsy \scriptfont\tw@=\sixsy \scriptscriptfont\tw@=\fivesy
  \textfont\thr@@=\eightex \scriptfont\thr@@=\sevenex
   \scriptscriptfont\thr@@=\sevenex
  \textfont\itfam=\eightit \scriptfont\itfam=\sevenit
   \scriptscriptfont\itfam=\sevenit
  \textfont\bffam=\eightbf \scriptfont\bffam=\sixbf
   \scriptscriptfont\bffam=\fivebf
 \setbox\strutbox\hbox{\vrule height7\p@ depth3\p@ width\z@}%
 \setbox\strutbox@\hbox{\raise.5\normallineskiplimit\vbox{%
   \kern-\normallineskiplimit\copy\strutbox}}%
 \setbox\z@\vbox{\hbox{$($}\kern\z@}\bigsize@=1.2\ht\z@
 \fi
 \normalbaselines\eightrm\ex@.2326ex\jot3\ex@\the\eightpoint@}

\font@\twelverm=cmr10 scaled\magstep1
\font@\twelveit=cmti10 scaled\magstep1
\font@\twelvesl=cmsl10 scaled\magstep1
\font@\twelvesmc=cmcsc10 scaled\magstep1
\font@\twelvett=cmtt10 scaled\magstep1
\font@\twelvebf=cmbx10 scaled\magstep1
\font@\twelvei=cmmi10 scaled\magstep1
\font@\twelvesy=cmsy10 scaled\magstep1
\font@\twelveex=cmex10 scaled\magstep1
\font@\twelvemsa=msam10 scaled\magstep1
\font@\twelveeufm=eufm10 scaled\magstep1
\font@\twelvemsb=msbm10 scaled\magstep1
\newtoks\twelvepoint@
\def\twelvepoint{\normalbaselineskip15\p@
 \abovedisplayskip15\p@ plus3.6\p@ minus10.8\p@
 \belowdisplayskip\abovedisplayskip
 \abovedisplayshortskip\z@ plus3.6\p@
 \belowdisplayshortskip8.4\p@ plus3.6\p@ minus4.8\p@
 \textonlyfont@\rm\twelverm \textonlyfont@\it\twelveit
 \textonlyfont@\sl\twelvesl \textonlyfont@\bf\twelvebf
 \textonlyfont@\smc\twelvesmc \textonlyfont@\tt\twelvett
 \textonlyfont@\bsmc\twelvebsmc
 \ifsyntax@ \def\big##1{{\hbox{$\left##1\right.$}}}%
  \let\Big\big \let\bigg\big \let\Bigg\big
 \else
  \textfont\z@=\twelverm  \scriptfont\z@=\tenrm  \scriptscriptfont\z@=\sevenrm
  \textfont\@ne=\twelvei  \scriptfont\@ne=\teni  \scriptscriptfont\@ne=\seveni
  \textfont\tw@=\twelvesy \scriptfont\tw@=\tensy \scriptscriptfont\tw@=\sevensy
  \textfont\thr@@=\twelveex \scriptfont\thr@@=\tenex
        \scriptscriptfont\thr@@=\tenex
  \textfont\itfam=\twelveit \scriptfont\itfam=\tenit
        \scriptscriptfont\itfam=\tenit
  \textfont\bffam=\twelvebf \scriptfont\bffam=\tenbf
        \scriptscriptfont\bffam=\sevenbf
  \setbox\strutbox\hbox{\vrule height10.2\p@ depth4.2\p@ width\z@}%
  \setbox\strutbox@\hbox{\lower.6\normallineskiplimit\vbox{%
        \kern-\normallineskiplimit\copy\strutbox}}%
 \setbox\z@\vbox{\hbox{$($}\kern\z@}\bigsize@=1.4\ht\z@
 \fi
 \normalbaselines\rm\ex@.2326ex\jot3.6\ex@\the\twelvepoint@}

\def\headfonts{\twelvepoint\bf}

\font@\fourteenrm=cmr10 scaled\magstep2
\font@\fourteenit=cmti10 scaled\magstep2
\font@\fourteensl=cmsl10 scaled\magstep2
\font@\fourteensmc=cmcsc10 scaled\magstep2
\font@\fourteentt=cmtt10 scaled\magstep2
\font@\fourteenbf=cmbx10 scaled\magstep2
\font@\fourteeni=cmmi10 scaled\magstep2
\font@\fourteensy=cmsy10 scaled\magstep2
\font@\fourteenex=cmex10 scaled\magstep2
\font@\fourteenmsa=msam10 scaled\magstep2
\font@\fourteeneufm=eufm10 scaled\magstep2
\font@\fourteenmsb=msbm10 scaled\magstep2
\newtoks\fourteenpoint@
\def\fourteenpoint{\normalbaselineskip15\p@
 \abovedisplayskip18\p@ plus4.3\p@ minus12.9\p@
 \belowdisplayskip\abovedisplayskip
 \abovedisplayshortskip\z@ plus4.3\p@
 \belowdisplayshortskip10.1\p@ plus4.3\p@ minus5.8\p@
 \textonlyfont@\rm\fourteenrm \textonlyfont@\it\fourteenit
 \textonlyfont@\sl\fourteensl \textonlyfont@\bf\fourteenbf
 \textonlyfont@\smc\fourteensmc \textonlyfont@\tt\fourteentt
 \textonlyfont@\bsmc\fourteenbsmc
 \ifsyntax@ \def\big##1{{\hbox{$\left##1\right.$}}}%
  \let\Big\big \let\bigg\big \let\Bigg\big
 \else
  \textfont\z@=\fourteenrm  \scriptfont\z@=\twelverm  \scriptscriptfont\z@=\tenrm
  \textfont\@ne=\fourteeni  \scriptfont\@ne=\twelvei  \scriptscriptfont\@ne=\teni
  \textfont\tw@=\fourteensy \scriptfont\tw@=\twelvesy \scriptscriptfont\tw@=\tensy
  \textfont\thr@@=\fourteenex \scriptfont\thr@@=\twelveex
        \scriptscriptfont\thr@@=\twelveex
  \textfont\itfam=\fourteenit \scriptfont\itfam=\twelveit
        \scriptscriptfont\itfam=\twelveit
  \textfont\bffam=\fourteenbf \scriptfont\bffam=\twelvebf
        \scriptscriptfont\bffam=\tenbf
  \setbox\strutbox\hbox{\vrule height12.2\p@ depth5\p@ width\z@}%
  \setbox\strutbox@\hbox{\lower.72\normallineskiplimit\vbox{%
        \kern-\normallineskiplimit\copy\strutbox}}%
 \setbox\z@\vbox{\hbox{$($}\kern\z@}\bigsize@=1.7\ht\z@
 \fi
 \normalbaselines\rm\ex@.2326ex\jot4.3\ex@\the\fourteenpoint@}

\def\chapheadfonts{\fourteenpoint\bf}

\font@\seventeenrm=cmr10 scaled\magstep3
\font@\seventeenit=cmti10 scaled\magstep3
\font@\seventeensl=cmsl10 scaled\magstep3
\font@\seventeensmc=cmcsc10 scaled\magstep3
\font@\seventeentt=cmtt10 scaled\magstep3
\font@\seventeenbf=cmbx10 scaled\magstep3
\font@\seventeeni=cmmi10 scaled\magstep3
\font@\seventeensy=cmsy10 scaled\magstep3
\font@\seventeenex=cmex10 scaled\magstep3
\font@\seventeenmsa=msam10 scaled\magstep3
\font@\seventeeneufm=eufm10 scaled\magstep3
\font@\seventeenmsb=msbm10 scaled\magstep3
\newtoks\seventeenpoint@
\def\seventeenpoint{\normalbaselineskip18\p@
 \abovedisplayskip21.6\p@ plus5.2\p@ minus15.4\p@
 \belowdisplayskip\abovedisplayskip
 \abovedisplayshortskip\z@ plus5.2\p@
 \belowdisplayshortskip12.1\p@ plus5.2\p@ minus7\p@
 \textonlyfont@\rm\seventeenrm \textonlyfont@\it\seventeenit
 \textonlyfont@\sl\seventeensl \textonlyfont@\bf\seventeenbf
 \textonlyfont@\smc\seventeensmc \textonlyfont@\tt\seventeentt
 \textonlyfont@\bsmc\seventeenbsmc
 \ifsyntax@ \def\big##1{{\hbox{$\left##1\right.$}}}%
  \let\Big\big \let\bigg\big \let\Bigg\big
 \else
  \textfont\z@=\seventeenrm  \scriptfont\z@=\fourteenrm  \scriptscriptfont\z@=\twelverm
  \textfont\@ne=\seventeeni  \scriptfont\@ne=\fourteeni  \scriptscriptfont\@ne=\twelvei
  \textfont\tw@=\seventeensy \scriptfont\tw@=\fourteensy \scriptscriptfont\tw@=\twelvesy
  \textfont\thr@@=\seventeenex \scriptfont\thr@@=\fourteenex
        \scriptscriptfont\thr@@=\fourteenex
  \textfont\itfam=\seventeenit \scriptfont\itfam=\fourteenit
        \scriptscriptfont\itfam=\fourteenit
  \textfont\bffam=\seventeenbf \scriptfont\bffam=\fourteenbf
        \scriptscriptfont\bffam=\twelvebf
  \setbox\strutbox\hbox{\vrule height14.6\p@ depth6\p@ width\z@}%
  \setbox\strutbox@\hbox{\lower.86\normallineskiplimit\vbox{%
        \kern-\normallineskiplimit\copy\strutbox}}%
 \setbox\z@\vbox{\hbox{$($}\kern\z@}\bigsize@=2\ht\z@
 \fi
 \normalbaselines\rm\ex@.2326ex\jot5.2\ex@\the\seventeenpoint@}

\font@\rrrrrm=cmr10 scaled\magstep4
\font@\bigtitlefont=cmbx10 scaled\magstep4

\parindent1pc
\normallineskiplimit\p@
\newdimen\indenti \indenti=2pc
\def\pageheight#1{\vsize#1}
\def\pagewidth#1{\hsize#1%
   \captionwidth@\hsize \advance\captionwidth@-2\indenti}
\pagewidth{30pc} \pageheight{47pc}
\def\topmatter{%
 \ifx\undefined\msafam
 \else\font@\eightmsa=msam8 \font@\sixmsa=msam6
   \ifsyntax@\else \addto\tenpoint{\textfont\msafam=\tenmsa
              \scriptfont\msafam=\sevenmsa \scriptscriptfont\msafam=\fivemsa}%
     \addto\eightpoint{\textfont\msafam=\eightmsa \scriptfont\msafam=\sixmsa
              \scriptscriptfont\msafam=\fivemsa}%
   \fi
 \fi
 \ifx\undefined\msbfam
 \else\font@\eightmsb=msbm8 \font@\sixmsb=msbm6
   \ifsyntax@\else \addto\tenpoint{\textfont\msbfam=\tenmsb
         \scriptfont\msbfam=\sevenmsb \scriptscriptfont\msbfam=\fivemsb}%
     \addto\eightpoint{\textfont\msbfam=\eightmsb \scriptfont\msbfam=\sixmsb
         \scriptscriptfont\msbfam=\fivemsb}%
   \fi
 \fi
 \ifx\undefined\eufmfam
 \else \font@\eighteufm=eufm8 \font@\sixeufm=eufm6
   \ifsyntax@\else \addto\tenpoint{\textfont\eufmfam=\teneufm
       \scriptfont\eufmfam=\seveneufm \scriptscriptfont\eufmfam=\fiveeufm}%
     \addto\eightpoint{\textfont\eufmfam=\eighteufm
       \scriptfont\eufmfam=\sixeufm \scriptscriptfont\eufmfam=\fiveeufm}%
   \fi
 \fi
 \ifx\undefined\eufbfam
 \else \font@\eighteufb=eufb8 \font@\sixeufb=eufb6
   \ifsyntax@\else \addto\tenpoint{\textfont\eufbfam=\teneufb
      \scriptfont\eufbfam=\seveneufb \scriptscriptfont\eufbfam=\fiveeufb}%
    \addto\eightpoint{\textfont\eufbfam=\eighteufb
      \scriptfont\eufbfam=\sixeufb \scriptscriptfont\eufbfam=\fiveeufb}%
   \fi
 \fi
 \ifx\undefined\eusmfam
 \else \font@\eighteusm=eusm8 \font@\sixeusm=eusm6
   \ifsyntax@\else \addto\tenpoint{\textfont\eusmfam=\teneusm
       \scriptfont\eusmfam=\seveneusm \scriptscriptfont\eusmfam=\fiveeusm}%
     \addto\eightpoint{\textfont\eusmfam=\eighteusm
       \scriptfont\eusmfam=\sixeusm \scriptscriptfont\eusmfam=\fiveeusm}%
   \fi
 \fi
 \ifx\undefined\eusbfam
 \else \font@\eighteusb=eusb8 \font@\sixeusb=eusb6
   \ifsyntax@\else \addto\tenpoint{\textfont\eusbfam=\teneusb
       \scriptfont\eusbfam=\seveneusb \scriptscriptfont\eusbfam=\fiveeusb}%
     \addto\eightpoint{\textfont\eusbfam=\eighteusb
       \scriptfont\eusbfam=\sixeusb \scriptscriptfont\eusbfam=\fiveeusb}%
   \fi
 \fi
 \ifx\undefined\eurmfam
 \else \font@\eighteurm=eurm8 \font@\sixeurm=eurm6
   \ifsyntax@\else \addto\tenpoint{\textfont\eurmfam=\teneurm
       \scriptfont\eurmfam=\seveneurm \scriptscriptfont\eurmfam=\fiveeurm}%
     \addto\eightpoint{\textfont\eurmfam=\eighteurm
       \scriptfont\eurmfam=\sixeurm \scriptscriptfont\eurmfam=\fiveeurm}%
   \fi
 \fi
 \ifx\undefined\eurbfam
 \else \font@\eighteurb=eurb8 \font@\sixeurb=eurb6
   \ifsyntax@\else \addto\tenpoint{\textfont\eurbfam=\teneurb
       \scriptfont\eurbfam=\seveneurb \scriptscriptfont\eurbfam=\fiveeurb}%
    \addto\eightpoint{\textfont\eurbfam=\eighteurb
       \scriptfont\eurbfam=\sixeurb \scriptscriptfont\eurbfam=\fiveeurb}%
   \fi
 \fi
 \ifx\undefined\cmmibfam
 \else \font@\eightcmmib=cmmib8 \font@\sixcmmib=cmmib6
   \ifsyntax@\else \addto\tenpoint{\textfont\cmmibfam=\tencmmib
       \scriptfont\cmmibfam=\sevencmmib \scriptscriptfont\cmmibfam=\fivecmmib}%
    \addto\eightpoint{\textfont\cmmibfam=\eightcmmib
       \scriptfont\cmmibfam=\sixcmmib \scriptscriptfont\cmmibfam=\fivecmmib}%
   \fi
 \fi
 \ifx\undefined\cmbsyfam
 \else \font@\eightcmbsy=cmbsy8 \font@\sixcmbsy=cmbsy6
   \ifsyntax@\else \addto\tenpoint{\textfont\cmbsyfam=\tencmbsy
      \scriptfont\cmbsyfam=\sevencmbsy \scriptscriptfont\cmbsyfam=\fivecmbsy}%
    \addto\eightpoint{\textfont\cmbsyfam=\eightcmbsy
      \scriptfont\cmbsyfam=\sixcmbsy \scriptscriptfont\cmbsyfam=\fivecmbsy}%
   \fi
 \fi
 \let\topmatter\relax}
\def\chapterno@{\uppercase\expandafter{\romannumeral\chaptercount@}}
\newcount\chaptercount@
\def\chapter{\nofrills@{\afterassignment\chapterno@
                        CHAPTER \global\chaptercount@=}\chapter@
 \DNii@##1{\leavevmode\hskip-\leftskip
   \rlap{\vbox to\z@{\vss\centerline{\eightpoint
   \chapter@##1\unskip}\baselineskip2pc\null}}\hskip\leftskip
   \nofrills@false}%
 \FN@\next@}
\newbox\titlebox@

\def\title{\nofrills@{\relax}\title@%
 \DNii@##1\endtitle{\global\setbox\titlebox@\vtop{\tenpoint\bf
 \raggedcenter@\ignorespaces
 \baselineskip1.3\baselineskip\title@{##1}\endgraf}%
 \ifmonograph@ \edef\next{\the\leftheadtoks}\ifx\next\empty
    \leftheadtext{##1}\fi
 \fi
 \edef\next{\the\rightheadtoks}\ifx\next\empty \rightheadtext{##1}\fi
 }\FN@\next@}
\newbox\authorbox@
\def\author#1\endauthor{\global\setbox\authorbox@
 \vbox{\tenpoint\smc\raggedcenter@\ignorespaces
 #1\endgraf}\relaxnext@ \edef\next{\the\leftheadtoks}%
 \ifx\next\empty\leftheadtext{#1}\fi}
\newbox\affilbox@
\def\affil#1\endaffil{\global\setbox\affilbox@
 \vbox{\tenpoint\raggedcenter@\ignorespaces#1\endgraf}}
\newcount\addresscount@
\addresscount@\z@
\def\address#1\endaddress{\global\advance\addresscount@\@ne
  \expandafter\gdef\csname address\number\addresscount@\endcsname
  {\vskip12\p@ minus6\p@\noindent\eightpoint\smc\ignorespaces#1\par}}
\def\email{\nofrills@{\eightpoint{\it E-mail\/}:\enspace}\email@
  \DNii@##1\endemail{%
  \expandafter\gdef\csname email\number\addresscount@\endcsname
  {\def\usualspace{{\it\enspace}}\smallskip\noindent\eightpoint\email@
  \ignorespaces##1\par}}%
 \FN@\next@}
\def\thedate@{}
\def\date#1\enddate{\gdef\thedate@{\tenpoint\ignorespaces#1\unskip}}
\def\thethanks@{}
\def\thanks#1\endthanks{\gdef\thethanks@{\eightpoint\ignorespaces#1.\unskip}}
\def\thekeywords@{}
\def\keywords{\nofrills@{{\it Key words and phrases.\enspace}}\keywords@
 \DNii@##1\endkeywords{\def\thekeywords@{\def\usualspace{{\it\enspace}}%
 \eightpoint\keywords@\ignorespaces##1\unskip.}}%
 \FN@\next@}
\def\thesubjclass@{}
\def\subjclass{\nofrills@{{\rm2000 {\it Mathematics Subject
   Classification\/}.\enspace}}\subjclass@
 \DNii@##1\endsubjclass{\def\thesubjclass@{\def\usualspace
  {{\rm\enspace}}\eightpoint\subjclass@\ignorespaces##1\unskip.}}%
 \FN@\next@}
\newbox\abstractbox@
\def\abstract{\nofrills@{{\smc Abstract.\enspace}}\abstract@
 \DNii@{\setbox\abstractbox@\vbox\bgroup\noindent$$\vbox\bgroup
  \def\envir@{abstract}\advance\hsize-2\indenti
  \usualspace@{{\enspace}}\eightpoint \noindent\abstract@\ignorespaces}%
 \FN@\next@}
\def\endabstract{\par\unskip\egroup$$\egroup}
\def\widestnumber#1#2{\begingroup\let\head\null\let\subhead\empty
   \let\subsubhead\subhead
   \ifx#1\head\global\setbox\tocheadbox@\hbox{#2.\enspace}%
   \else\ifx#1\subhead\global\setbox\tocsubheadbox@\hbox{#2.\enspace}%
   \else\ifx#1\key\bgroup\let\endrefitem@\egroup
        \key#2\endrefitem@\global\refindentwd\wd\keybox@
   \else\ifx#1\no\bgroup\let\endrefitem@\egroup
        \no#2\endrefitem@\global\refindentwd\wd\nobox@
   \else\ifx#1\page\global\setbox\pagesbox@\hbox{\quad\bf#2}%
   \else\ifx#1\item\setboxz@h{#2}\global\rosteritemwd\wdz@
        \global\advance\rosteritemwd by.5\parindent
   \else\message{\string\widestnumber is not defined for this option
   (\string#1)}%
\fi\fi\fi\fi\fi\fi\endgroup}
\newif\ifmonograph@
\def\Monograph{\monograph@true \let\headmark\rightheadtext
  \let\varindent@\indent \def\headfont@{\bf}\def\proclaimheadfont@{\smc}%
  \def\demofont@{\smc}}
\let\varindent@\indent

\newbox\tocheadbox@    \newbox\tocsubheadbox@
\newbox\tocbox@
\def\toc{\toc@{Contents}}
\def\newtocdefs{%
   \def \title##1\endtitle
       {\penaltyandskip@\z@\smallskipamount
        \hangindent\wd\tocheadbox@\noindent{\bf##1}}%
   \def \chapter##1{%
        Chapter \uppercase\expandafter{\romannumeral##1.\unskip}\enspace}%
   \def \specialhead##1\endspecialhead
       {\par\hangindent\wd\tocheadbox@ \noindent##1\par}%
   \def \head##1 ##2\endhead
       {\par\hangindent\wd\tocheadbox@ \noindent
        \if\notempty{##1}\hbox to\wd\tocheadbox@{\hfil##1\enspace}\fi
        ##2\par}%
   \def \subhead##1 ##2\endsubhead
       {\par\vskip-\parskip {\normalbaselines
        \advance\leftskip\wd\tocheadbox@
        \hangindent\wd\tocsubheadbox@ \noindent
        \if\notempty{##1}\hbox to\wd\tocsubheadbox@{##1\unskip\hfil}\fi
         ##2\par}}%
   \def \subsubhead##1 ##2\endsubsubhead
       {\par\vskip-\parskip {\normalbaselines
        \advance\leftskip\wd\tocheadbox@
        \hangindent\wd\tocsubheadbox@ \noindent
        \if\notempty{##1}\hbox to\wd\tocsubheadbox@{##1\unskip\hfil}\fi
        ##2\par}}}
\def\toc@#1{\relaxnext@
   \def\page##1%
       {\unskip\penalty0\null\hfil
        \rlap{\hbox to\wd\pagesbox@{\quad\hfil##1}}\hfilneg\penalty\@M}%
 \DN@{\ifx\next\nofrills\DN@\nofrills{\nextii@}%
      \else\DN@{\nextii@{{#1}}}\fi
      \next@}%
 \DNii@##1{%
\ifmonograph@\bgroup\else\setbox\tocbox@\vbox\bgroup
   \centerline{\headfont@\ignorespaces##1\unskip}\nobreak
   \vskip\belowheadskip \fi
   \setbox\tocheadbox@\hbox{0.\enspace}%
   \setbox\tocsubheadbox@\hbox{0.0.\enspace}%
   \leftskip\indenti \rightskip\leftskip
   \setbox\pagesbox@\hbox{\bf\quad000}\advance\rightskip\wd\pagesbox@
   \newtocdefs
 }%
 \FN@\next@}
\def\endtoc{\par\egroup}
\let\pretitle\relax
\let\preauthor\relax
\let\preaffil\relax
\let\predate\relax
\let\preabstract\relax
\let\prepaper\relax
\def\dedicatory #1\enddedicatory{\def\preabstract{{\medskip
  \eightpoint\it \raggedcenter@#1\endgraf}}}
\def\thetranslator@{}
\def\translator#1\endtranslator{\def\thetranslator@{\nobreak\medskip
 \line{\eightpoint\hfil Translated by \uppercase{#1}\qquad\qquad}\nobreak}}
\outer\def\endtopmatter{\runaway@{abstract}%
 \edef\next{\the\leftheadtoks}\ifx\next\empty
  \expandafter\leftheadtext\expandafter{\the\rightheadtoks}\fi
 \ifmonograph@\else
   \ifx\thesubjclass@\empty\else \makefootnote@{}{\thesubjclass@}\fi
   \ifx\thekeywords@\empty\else \makefootnote@{}{\thekeywords@}\fi
   \ifx\thethanks@\empty\else \makefootnote@{}{\thethanks@}\fi
 \fi
  \pretitle
  \ifmonograph@ \topskip7pc \else \topskip4pc \fi
  \box\titlebox@
  \topskip10pt
  \preauthor
  \ifvoid\authorbox@\else \vskip2.5pc plus1pc \unvbox\authorbox@\fi
  \preaffil
  \ifvoid\affilbox@\else \vskip1pc plus.5pc \unvbox\affilbox@\fi
  \predate
  \ifx\thedate@\empty\else \vskip1pc plus.5pc \line{\hfil\thedate@\hfil}\fi
  \preabstract
  \ifvoid\abstractbox@\else \vskip1.5pc plus.5pc \unvbox\abstractbox@ \fi
  \ifvoid\tocbox@\else\vskip1.5pc plus.5pc \unvbox\tocbox@\fi
  \prepaper
  \vskip2pc plus1pc
}
\def\document{\let\fontlist@\relax\let\alloclist@\relax
  \tenpoint}

\newskip\aboveheadskip       \aboveheadskip1.8\bigskipamount
\newdimen\belowheadskip      \belowheadskip1.8\medskipamount

\def\headfont@{\smc}
\def\penaltyandskip@#1#2{\relax\ifdim\lastskip<#2\relax\removelastskip
      \ifnum#1=\z@\else\penalty@#1\relax\fi\vskip#2%
  \else\ifnum#1=\z@\else\penalty@#1\relax\fi\fi}
\def\nobreak{\penalty\@M
  \ifvmode\def\penalty@{\let\penalty@\penalty\count@@@}%
  \everypar{\let\penalty@\penalty\everypar{}}\fi}
\let\penalty@\penalty
\def\heading#1\endheading{\head#1\endhead}

\def\specialheadfont@{\bf}
\outer\def\specialhead{\par\penaltyandskip@{-200}\aboveheadskip
  \begingroup\interlinepenalty\@M\rightskip\z@ plus\hsize \let\\\linebreak
  \specialheadfont@\noindent\ignorespaces}
\def\endspecialhead{\par\endgroup\nobreak\vskip\belowheadskip}
\let\headmark\eat@
\newskip\subheadskip       \subheadskip\medskipamount
\def\subheadfont@{\bf}
\outer\def\subhead{\nofrills@{.\enspace}\subhead@
 \DNii@##1\endsubhead{\par\penaltyandskip@{-100}\subheadskip
  \varindent@{\usualspace@{{\subheadfont@\enspace}}%
 \subheadfont@\ignorespaces##1\unskip\subhead@}\ignorespaces}%
 \FN@\next@}
\outer\def\subsubhead{\nofrills@{.\enspace}\subsubhead@
 \DNii@##1\endsubsubhead{\par\penaltyandskip@{-50}\medskipamount
      {\usualspace@{{\it\enspace}}%
  \it\ignorespaces##1\unskip\subsubhead@}\ignorespaces}%
 \FN@\next@}
\def\proclaimheadfont@{\bf}
\outer\def\proclaim{\runaway@{proclaim}\def\envir@{proclaim}%
  \nofrills@{.\enspace}\proclaim@
 \DNii@##1{\penaltyandskip@{-100}\medskipamount\varindent@
   \usualspace@{{\proclaimheadfont@\enspace}}\proclaimheadfont@
   \ignorespaces##1\unskip\proclaim@
  \sl\ignorespaces}%
 \FN@\next@}
\outer\def\endproclaim{\let\envir@\relax\par\rm
  \penaltyandskip@{55}\medskipamount}
\def\demoheadfont@{\it}
\def\demo{\runaway@{proclaim}\nofrills@{.\enspace}\demo@
     \DNii@##1{\par\penaltyandskip@\z@\medskipamount
  {\usualspace@{{\demoheadfont@\enspace}}%
  \varindent@\demoheadfont@\ignorespaces##1\unskip\demo@}\rm
  \ignorespaces}\FN@\next@}
\def\enddemo{\par\medskip}
\def\qed{\ifhmode\unskip\nobreak\fi\quad\ifmmode\square\else$\m@th\square$\fi}
\let\remark\demo
\let\endremark\enddemo
\def\definition{\runaway@{proclaim}%
  \nofrills@{.\demoheadfont@\enspace}\definition@
        \DNii@##1{\penaltyandskip@{-100}\medskipamount
        {\usualspace@{{\demoheadfont@\enspace}}%
        \varindent@\demoheadfont@\ignorespaces##1\unskip\definition@}%
        \rm \ignorespaces}\FN@\next@}


\newdimen\rosteritemwd
\newcount\rostercount@
\newif\iffirstitem@
\let\plainitem@\item
\newtoks\everypartoks@
\def\par@{\everypartoks@\expandafter{\the\everypar}\everypar{}}
\def\roster{\edef\leftskip@{\leftskip\the\leftskip}%
 \relaxnext@
 \rostercount@\z@  
 \def\item{\FN@\rosteritem@}%
 \DN@{\ifx\next\runinitem\let\next@\nextii@\else
  \let\next@\nextiii@\fi\next@}%
 \DNii@\runinitem  
  {\unskip  
   \DN@{\ifx\next[\let\next@\nextii@\else
    \ifx\next"\let\next@\nextiii@\else\let\next@\nextiv@\fi\fi\next@}%
   \DNii@[####1]{\rostercount@####1\relax
    \enspace{\rm(\number\rostercount@)}~\ignorespaces}%
   \def\nextiii@"####1"{\enspace{\rm####1}~\ignorespaces}%
   \def\nextiv@{\enspace{\rm(1)}\rostercount@\@ne~}%
   \par@\firstitem@false  
   \FN@\next@}%
 \def\nextiii@{\par\par@  
  \penalty\@m\smallskip\vskip-\parskip
  \firstitem@true}%
 \FN@\next@}
\def\rosteritem@{\iffirstitem@\firstitem@false\else\par\vskip-\parskip\fi
 \leftskip3\parindent\noindent  
 \DNii@[##1]{\rostercount@##1\relax
  \llap{\hbox to2.5\parindent{\hss\rm(\number\rostercount@)}%
   \hskip.5\parindent}\ignorespaces}%
 \def\nextiii@"##1"{%
  \llap{\hbox to2.5\parindent{\hss\rm##1}\hskip.5\parindent}\ignorespaces}%
 \def\nextiv@{\advance\rostercount@\@ne
  \llap{\hbox to2.5\parindent{\hss\rm(\number\rostercount@)}%
   \hskip.5\parindent}}%
 \ifx\next[\let\next@\nextii@\else\ifx\next"\let\next@\nextiii@\else
  \let\next@\nextiv@\fi\fi\next@}

\newif\ifnextRunin@
\def\endroster{\relaxnext@
 \par\leftskip@  
 \penalty-50 \vskip-\parskip\smallskip  
 \DN@{\ifx\next\Runinitem\let\next@\relax
  \else\nextRunin@false\let\item\plainitem@  
   \ifx\next\par 
    \DN@\par{\everypar\expandafter{\the\everypartoks@}}%
   \else  
    \DN@{\noindent\everypar\expandafter{\the\everypartoks@}}%
  \fi\fi\next@}%
 \FN@\next@}
\newcount\rosterhangafter@
\def\Runinitem#1\roster\runinitem{\relaxnext@
 \rostercount@\z@ 
 \def\item{\FN@\rosteritem@}%
 \def\runinitem@{#1}%
 \DN@{\ifx\next[\let\next\nextii@\else\ifx\next"\let\next\nextiii@
  \else\let\next\nextiv@\fi\fi\next}%
 \DNii@[##1]{\rostercount@##1\relax
  \def\item@{{\rm(\number\rostercount@)}}\nextv@}%
 \def\nextiii@"##1"{\def\item@{{\rm##1}}\nextv@}%
 \def\nextiv@{\advance\rostercount@\@ne
  \def\item@{{\rm(\number\rostercount@)}}\nextv@}%
 \def\nextv@{\setbox\z@\vbox  
  {\ifnextRunin@\noindent\fi  
  \runinitem@\unskip\enspace\item@~\par  
  \global\rosterhangafter@\prevgraf}%
  \firstitem@false  
  \ifnextRunin@\else\par\fi  
  \hangafter\rosterhangafter@\hangindent3\parindent
  \ifnextRunin@\noindent\fi  
  \runinitem@\unskip\enspace 
  \item@~\ifnextRunin@\else\par@\fi  
  \nextRunin@true\ignorespaces}%
 \FN@\next@}
\def\footmarkform@#1{$\m@th^{#1}$}
\let\thefootnotemark\footmarkform@
\def\makefootnote@#1#2{\insert\footins
 {\interlinepenalty\interfootnotelinepenalty
 \eightpoint\splittopskip\ht\strutbox\splitmaxdepth\dp\strutbox
 \floatingpenalty\@MM\leftskip\z@\rightskip\z@\spaceskip\z@\xspaceskip\z@
 \leavevmode{#1}\footstrut\ignorespaces#2\unskip\lower\dp\strutbox
 \vbox to\dp\strutbox{}}}
\newcount\footmarkcount@
\footmarkcount@\z@
\def\footnotemark{\let\@sf\empty\relaxnext@
 \ifhmode\edef\@sf{\spacefactor\the\spacefactor}\/\fi
 \DN@{\ifx[\next\let\next@\nextii@\else
  \ifx"\next\let\next@\nextiii@\else
  \let\next@\nextiv@\fi\fi\next@}%
 \DNii@[##1]{\footmarkform@{##1}\@sf}%
 \def\nextiii@"##1"{{##1}\@sf}%
 \def\nextiv@{\iffirstchoice@\global\advance\footmarkcount@\@ne\fi
  \footmarkform@{\number\footmarkcount@}\@sf}%
 \FN@\next@}
\def\footnotetext{\relaxnext@
 \DN@{\ifx[\next\let\next@\nextii@\else
  \ifx"\next\let\next@\nextiii@\else
  \let\next@\nextiv@\fi\fi\next@}%
 \DNii@[##1]##2{\makefootnote@{\footmarkform@{##1}}{##2}}%
 \def\nextiii@"##1"##2{\makefootnote@{##1}{##2}}%
 \def\nextiv@##1{\makefootnote@{\footmarkform@{\number\footmarkcount@}}{##1}}%
 \FN@\next@}
\def\footnote{\let\@sf\empty\relaxnext@
 \ifhmode\edef\@sf{\spacefactor\the\spacefactor}\/\fi
 \DN@{\ifx[\next\let\next@\nextii@\else
  \ifx"\next\let\next@\nextiii@\else
  \let\next@\nextiv@\fi\fi\next@}%
 \DNii@[##1]##2{\footnotemark[##1]\footnotetext[##1]{##2}}%
 \def\nextiii@"##1"##2{\footnotemark"##1"\footnotetext"##1"{##2}}%
 \def\nextiv@##1{\footnotemark\footnotetext{##1}}%
 \FN@\next@}
\def\adjustfootnotemark#1{\advance\footmarkcount@#1\relax}
\def\footnoterule{\kern-3\p@
  \hrule width 5pc\kern 2.6\p@} 
\def\captionfont@{\smc}
\def\topcaption#1#2\endcaption{%
  {\dimen@\hsize \advance\dimen@-\captionwidth@
   \rm\raggedcenter@ \advance\leftskip.5\dimen@ \rightskip\leftskip
  {\captionfont@#1}%
  \if\notempty{#2}.\enspace\ignorespaces#2\fi
  \endgraf}\nobreak\bigskip}
\def\botcaption#1#2\endcaption{%
  \nobreak\bigskip
  \setboxz@h{\captionfont@#1\if\notempty{#2}.\enspace\rm#2\fi}%
  {\dimen@\hsize \advance\dimen@-\captionwidth@
   \leftskip.5\dimen@ \rightskip\leftskip
   \noindent \ifdim\wdz@>\captionwidth@ 
   \else\hfil\fi 
  {\captionfont@#1}\if\notempty{#2}.\enspace\rm#2\fi\endgraf}}
\def\@ins{\par\begingroup\def\vspace##1{\vskip##1\relax}%
  \def\captionwidth##1{\captionwidth@##1\relax}%
  \setbox\z@\vbox\bgroup} 
\def\block{\RIfMIfI@\nondmatherr@\block\fi
       \else\ifvmode\vskip\abovedisplayskip\noindent\fi
        $$\def\endblock{\par\egroup$$}\fi
  \vbox\bgroup\advance\hsize-2\indenti\noindent}
\def\endblock{\par\egroup}
\def\cite#1{{\rm[{\citefont@\m@th#1}]}}
\def\citefont@{\rm}
\def\refsfont@{\eightpoint}
\outer\def\Refs{\runaway@{proclaim}%
 \relaxnext@ \DN@{\ifx\next\nofrills\DN@\nofrills{\nextii@}\else
  \DN@{\nextii@{References}}\fi\next@}%
 \DNii@##1{\penaltyandskip@{-200}\aboveheadskip
  \line{\hfil\headfont@\ignorespaces##1\unskip\hfil}\nobreak
  \vskip\belowheadskip
  \begingroup\refsfont@\sfcode`.=\@m}%
 \FN@\next@}
\def\endRefs{\par\endgroup}
\newbox\nobox@            \newbox\keybox@           \newbox\bybox@
\newbox\paperbox@         \newbox\paperinfobox@     \newbox\jourbox@
\newbox\volbox@           \newbox\issuebox@         \newbox\yrbox@
\newbox\pagesbox@         \newbox\bookbox@          \newbox\bookinfobox@
\newbox\publbox@          \newbox\publaddrbox@      \newbox\finalinfobox@
\newbox\edsbox@           \newbox\langbox@
\newif\iffirstref@        \newif\iflastref@
\newif\ifprevjour@        \newif\ifbook@            \newif\ifprevinbook@
\newif\ifquotes@          \newif\ifbookquotes@      \newif\ifpaperquotes@
\newdimen\bysamerulewd@
\setboxz@h{\refsfont@\kern3em}
\bysamerulewd@\wdz@
\newdimen\refindentwd
\setboxz@h{\refsfont@ 00. }
\refindentwd\wdz@
\outer\def\ref{\begingroup \noindent\hangindent\refindentwd
 \firstref@true \def\nofrills{\def\refkern@{\kern3sp}}%
 \ref@}
\def\ref@{\book@false \bgroup\let\endrefitem@\egroup \ignorespaces}
\def\moreref{\endrefitem@\endref@\firstref@false\ref@}%
\def\transl{\endrefitem@\endref@\firstref@false
  \book@false
  \prepunct@
  \setboxz@h\bgroup \aftergroup\unhbox\aftergroup\z@
    \def\endrefitem@{\unskip\refkern@\egroup}\ignorespaces}%
\def\emptyifempty@{\dimen@\wd\currbox@
  \advance\dimen@-\wd\z@ \advance\dimen@-.1\p@
  \ifdim\dimen@<\z@ \setbox\currbox@\copy\voidb@x \fi}
\let\refkern@\relax
\def\endrefitem@{\unskip\refkern@\egroup
  \setboxz@h{\refkern@}\emptyifempty@}\ignorespaces
\def\refdef@#1#2#3{\edef\next@{\noexpand\endrefitem@
  \let\noexpand\currbox@\csname\expandafter\eat@\string#1box@\endcsname
    \noexpand\setbox\noexpand\currbox@\hbox\bgroup}%
  \toks@\expandafter{\next@}%
  \if\notempty{#2#3}\toks@\expandafter{\the\toks@
  \def\endrefitem@{\unskip#3\refkern@\egroup
  \setboxz@h{#2#3\refkern@}\emptyifempty@}#2}\fi
  \toks@\expandafter{\the\toks@\ignorespaces}%
  \edef#1{\the\toks@}}
\refdef@\no{}{. }
\refdef@\key{[\m@th}{] }
\refdef@\by{}{}
\def\bysame{\by\hbox to\bysamerulewd@{\hrulefill}\thinspace
   \kern0sp}
\def\manyby{\message{\string\manyby is no longer necessary; \string\by
  can be used instead, starting with version 2.0 of \styname.STY}\by}
\refdef@\paper{\ifpaperquotes@``\fi\it}{}
\refdef@\paperinfo{}{}
\def\jour{\endrefitem@\let\currbox@\jourbox@
  \setbox\currbox@\hbox\bgroup
  \def\endrefitem@{\unskip\refkern@\egroup
    \setboxz@h{\refkern@}\emptyifempty@
    \ifvoid\jourbox@\else\prevjour@true\fi}%
\ignorespaces}
\refdef@\vol{\ifbook@\else\bf\fi}{}
\refdef@\issue{no. }{}
\refdef@\yr{}{}
\refdef@\pages{}{}
\def\page{\endrefitem@\def\pp@{\def\pp@{pp.~}p.~}\let\currbox@\pagesbox@
  \setbox\currbox@\hbox\bgroup\ignorespaces}
\def\pp@{pp.~}
\def\book{\endrefitem@ \let\currbox@\bookbox@
 \setbox\currbox@\hbox\bgroup\def\endrefitem@{\unskip\refkern@\egroup
  \setboxz@h{\ifbookquotes@``\fi}\emptyifempty@
  \ifvoid\bookbox@\else\book@true\fi}%
  \ifbookquotes@``\fi\it\ignorespaces}
\def\inbook{\endrefitem@
  \let\currbox@\bookbox@\setbox\currbox@\hbox\bgroup
  \def\endrefitem@{\unskip\refkern@\egroup
  \setboxz@h{\ifbookquotes@``\fi}\emptyifempty@
  \ifvoid\bookbox@\else\book@true\previnbook@true\fi}%
  \ifbookquotes@``\fi\ignorespaces}
\refdef@\eds{(}{, eds.)}
\def\ed{\endrefitem@\let\currbox@\edsbox@
 \setbox\currbox@\hbox\bgroup
 \def\endrefitem@{\unskip, ed.)\refkern@\egroup
  \setboxz@h{(, ed.)}\emptyifempty@}(\ignorespaces}
\refdef@\bookinfo{}{}
\refdef@\publ{}{}
\refdef@\publaddr{}{}
\refdef@\finalinfo{}{}
\refdef@\lang{(}{)}
\def\toappear{\nofrills\finalinfo(to appear)}
\let\refdef@\relax 
\def\ppunbox@#1{\ifvoid#1\else\prepunct@\unhbox#1\fi}
\def\nocomma@#1{\ifvoid#1\else\changepunct@3\prepunct@\unhbox#1\fi}
\def\changepunct@#1{\ifnum\lastkern<3 \unkern\kern#1sp\fi}
\def\prepunct@{\count@\lastkern\unkern
  \ifnum\lastpenalty=0
    \let\penalty@\relax
  \else
    \edef\penalty@{\penalty\the\lastpenalty\relax}%
  \fi
  \unpenalty
  \let\refspace@\ \ifcase\count@,
\or;\or.\or 
  \or\let\refspace@\relax
  \else,\fi
  \ifquotes@''\quotes@false\fi \penalty@ \refspace@
}
\def\transferpenalty@#1{\dimen@\lastkern\unkern
  \ifnum\lastpenalty=0\unpenalty\let\penalty@\relax
  \else\edef\penalty@{\penalty\the\lastpenalty\relax}\unpenalty\fi
  #1\penalty@\kern\dimen@}
\def\endref{\endrefitem@\lastref@true\endref@
  \par\endgroup \prevjour@false \previnbook@false }
\def\endref@{%
\iffirstref@
  \ifvoid\nobox@\ifvoid\keybox@\indent\fi
  \else\hbox to\refindentwd{\hss\unhbox\nobox@}\fi
  \ifvoid\keybox@
  \else\ifdim\wd\keybox@>\refindentwd
         \box\keybox@
       \else\hbox to\refindentwd{\unhbox\keybox@\hfil}\fi\fi
  \kern4sp\ppunbox@\bybox@
\fi 
  \ifvoid\paperbox@
  \else\prepunct@\unhbox\paperbox@
    \ifpaperquotes@\quotes@true\fi\fi
  \ppunbox@\paperinfobox@
  \ifvoid\jourbox@
    \ifprevjour@ \nocomma@\volbox@
      \nocomma@\issuebox@
      \ifvoid\yrbox@\else\changepunct@3\prepunct@(\unhbox\yrbox@
        \transferpenalty@)\fi
      \ppunbox@\pagesbox@
    \fi 
  \else \prepunct@\unhbox\jourbox@
    \nocomma@\volbox@
    \nocomma@\issuebox@
    \ifvoid\yrbox@\else\changepunct@3\prepunct@(\unhbox\yrbox@
      \transferpenalty@)\fi
    \ppunbox@\pagesbox@
  \fi 
  \ifbook@\prepunct@\unhbox\bookbox@ \ifbookquotes@\quotes@true\fi \fi
  \nocomma@\edsbox@
  \ppunbox@\bookinfobox@
  \ifbook@\ifvoid\volbox@\else\prepunct@ vol.~\unhbox\volbox@
  \fi\fi
  \ppunbox@\publbox@ \ppunbox@\publaddrbox@
  \ifbook@ \ppunbox@\yrbox@
    \ifvoid\pagesbox@
    \else\prepunct@\pp@\unhbox\pagesbox@\fi
  \else
    \ifprevinbook@ \ppunbox@\yrbox@
      \ifvoid\pagesbox@\else\prepunct@\pp@\unhbox\pagesbox@\fi
    \fi \fi
  \ppunbox@\finalinfobox@
  \iflastref@
    \ifvoid\langbox@.\ifquotes@''\fi
    \else\changepunct@2\prepunct@\unhbox\langbox@\fi
  \else
    \ifvoid\langbox@\changepunct@1%
    \else\changepunct@3\prepunct@\unhbox\langbox@
      \changepunct@1\fi
  \fi
}
\outer\def\enddocument{%
 \runaway@{proclaim}%
\ifmonograph@ 
\else
 \nobreak
 \thetranslator@
 \count@\z@ \loop\ifnum\count@<\addresscount@\advance\count@\@ne
 \csname address\number\count@\endcsname
 \csname email\number\count@\endcsname
 \repeat
\fi
 \vfill\supereject\end}

\def\headfont@{\headfonts}
\def\proclaimheadfont@{\bf}
\def\specialheadfont@{\bf}
\def\subheadfont@{\bf}
\def\demoheadfont@{\smc}

\newif\ifThisToToc \ThisToTocfalse
\newif\iftocloaded \tocloadedfalse

\def\C@L{\noexpand\Cal}\def\B@B{\noexpand\Bbb}\def\fR@K{\noexpand\frak}
\def\S@{\noexpand\S}\def\P@P{\noexpand\"}
\def\xpar{\\}

\def\writetoc#1{\iftocloaded\ifThisToToc\begingroup\def\totoc{}
  \def\Cal{\noexpand\C@L}\def\Bbb{\noexpand\B@B}
  \def\frak{\noexpand\fR@K}\def\goth{\frak}\def\S{\noexpand\S@}
  \def\"{\noexpand\P@P}
  \def\xpar{\par\penalty100000 }\def\idx##1{##1}\def\\{\xpar}
  \edef\next@{\write\toc{\noindent#1\leaderfill\noexpand\folio\par}}%
  \next@\endgroup\global\ThisToTocfalse\fi\fi}
\def\leaderfill{\leaders\hbox to 1em{\hss.\hss}\hfill}

\newif\ifindexloaded \indexloadedfalse
\def\idx#1{\ifindexloaded\begingroup\def\ign{}\def\it{}\def\/{}%
 \def\smc{}\def\bf{}\def\tt{}%
 \def\Cal{\noexpand\C@L}\def\Bbb{\noexpand\B@B}%
 \def\frak{\noexpand\fR@K}\def\goth{\frak}\def\S{\noexpand\S@}%
  \def\"{\noexpand\P@P}%
 {\edef\next@{\write\index{#1, \noexpand\folio}}\next@}%
 \endgroup\fi{#1}}
\def\ign#1{}

\def\totoc{\global\ThisToToctrue}

\outer\def\head#1\endhead{\par\penaltyandskip@{-200}\aboveheadskip
 {\headfont@\raggedcenter@\interlinepenalty\@M
 \ignorespaces#1\endgraf}\nobreak
 \vskip\belowheadskip
 \headmark{#1}\writetoc{#1}}

\outer\def\chaphead#1\endchaphead{\par\penaltyandskip@{-200}\aboveheadskip
 {\chapheadfonts\raggedcenter@\interlinepenalty\@M
 \ignorespaces#1\endgraf}\nobreak
 \vskip3\belowheadskip
 \headmark{#1}\writetoc{#1}}

\def\folio{{\foliofont@\ifnum\pageno<\z@ \romannumeral-\pageno
 \else\number\pageno \fi}}
\newtoks\leftheadtoks
\newtoks\rightheadtoks

\def\leftheadtext{\nofrills@{\relax}\lht@
  \DNii@##1{\leftheadtoks\expandafter{\lht@{##1}}%
    \mark{\the\leftheadtoks\noexpand\else\the\rightheadtoks}
    \ifsyntax@\setboxz@h{\def\\{\unskip\space\ignorespaces}%
        \headlinefont@##1}\fi}%
  \FN@\next@}
\def\rightheadtext{\nofrills@{\relax}\rht@
  \DNii@##1{\rightheadtoks\expandafter{\rht@{##1}}%
    \mark{\the\leftheadtoks\noexpand\else\the\rightheadtoks}%
    \ifsyntax@\setboxz@h{\def\\{\unskip\space\ignorespaces}%
        \headlinefont@##1}\fi}%
  \FN@\next@}
\def\NoRunningHeads{\global\runheads@false\global\let\headmark\eat@}

\newif\iffirstpage@     \firstpage@true
\newif\ifrunheads@      \runheads@true

\newdimen\fullhsize \fullhsize=\hsize
\newdimen\fullvsize \fullvsize=\vsize
\def\fullline{\hbox to\fullhsize}

\def\pagenumbers{\gdef\folio{\folio@}}

\let\norunningheads\NoRunningHeads
\def\userunningheads{\global\runheads@true}
\norunningheads

\headline={\def\chapter#1{\chapterno@. }%
  \def\\{\unskip\space\ignorespaces}\ifrunheads@\headlinefont@
    \ifodd\pageno\rightheadline \else\leftheadline\fi
   \else\hfil\fi\ifNoRunHeadline\global\NoRunHeadlinefalse\fi}
\let\folio@\folio
\def\foliofont@{\foliofont}
\def\foliofont{\eightrm}
\def\headlinefont@{\headlinefont}
\def\headlinefont{\eightpoint\smc}
\def\leftheadline{\rlap{\folio}\hfill
   \ifNoRunHeadline\else\iftrue\topmark\fi\fi \hfill}
\def\rightheadline{\hfill\ifNoRunHeadline
   \else \expandafter\fi
  \hfill \llap{\folio}}
\footline={{\eightpoint\bottremark}%
   \ifrunheads@\else\hfil{\let\foliofont\tenrm\folio}\fi\hfil}
\def\bottremark{}
 
\newif\ifNoRunHeadline      
\def\norunninghead{\global\NoRunHeadlinetrue}
\norunninghead

\output={\output@}
%
\newif\ifoffset\offsetfalse
\output={\output@}
\def\output@{%
 \ifoffset 
  \ifodd\count0\advance\hoffset by0.5truecm
   \else\advance\hoffset by-0.5truecm\fi\fi
 \shipout\vbox{%
  \makeheadline \pagebody \makefootline }%
 \advancepageno \ifnum\outputpenalty>-\@MM\else\dosupereject\fi}

\def\indexoutput#1{%
  \ifoffset 
   \ifodd\count0\advance\hoffset by0.5truecm
    \else\advance\hoffset by-0.5truecm\fi\fi
  \shipout\vbox{\makeheadline
  \vbox to\fullvsize{\boxmaxdepth\maxdepth%
  \ifvoid\topins\else\unvbox\topins\fi%
  #1 %
  \ifvoid\footins\else 
    \vskip\skip\footins
    \footnoterule
    \unvbox\footins\fi
  \ifr@ggedbottom \kern-\dimen@ \vfil \fi}%
  \baselineskip2pc
  \makefootline}%
 \global\advance\pageno\@ne
 \ifnum\outputpenalty>-\@MM\else\dosupereject\fi}
 
 \newbox\partialpage \newdimen\halfsize \halfsize=0.5\fullhsize
 \advance\halfsize by-0.5em

 \def\begindoublecolumns{\output={\indexoutput{\unvbox255}}%
   \begingroup \def\line{\fullline}
   \output={\global\setbox\partialpage=\vbox{\unvbox255\bigskip}}\eject
   \output={\doublecolumnout}\hsize=\halfsize \vsize=2\fullvsize}
 \def\enddoublecolumns{\output={\balancecolumns}\eject
  \endgroup \pagegoal=\fullvsize%
  \output={\output@}}
\def\doublecolumnout{\splittopskip=\topskip \splitmaxdepth=\maxdepth
  \dimen@=\fullvsize \advance\dimen@ by-\ht\partialpage
  \setbox0=\vsplit255 to \dimen@ \setbox2=\vsplit255 to \dimen@
  \indexoutput{\pagesofar} \unvbox255 \penalty\outputpenalty}
\def\pagesofar{\unvbox\partialpage
  \wd0=\hsize \wd2=\hsize \hbox to\fullhsize{\box0\hfil\box2}}
\def\balancecolumns{\setbox0=\vbox{\unvbox255} \dimen@=\ht0
  \advance\dimen@ by\topskip \advance\dimen@ by-\baselineskip
  \divide\dimen@ by2 \splittopskip=\topskip
  {\vbadness=10000 \loop \global\setbox3=\copy0
    \global\setbox1=\vsplit3 to\dimen@
    \ifdim\ht3>\dimen@ \global\advance\dimen@ by1pt \repeat}
  \setbox0=\vbox to\dimen@{\unvbox1} \setbox2=\vbox to\dimen@{\unvbox3}
  \pagesofar}

\tenpoint
\catcode`\@=\active

\def\smallheadings{\let\chapheadfonts\tenpoint\let\headfonts\tenpoint}

\tenpoint
\catcode`\@=\active

\def\LL{\leavevmode\setbox0=\hbox{L}\hbox to\wd0{\hss\char'40L}}
\def\al{\alpha}
\def\be{\beta}


\def\today{\ifcase\month\or
 January\or February\or March\or April\or May\or June\or
 July\or August\or September\or October\or November\or December\fi
 \space\number\day, \number\year}

\def\({\left(}
\def\){\right)}
\def\[{\left[}
\def\]{\right]}

\def\3{\ss}
\catcode`\@=11
\def\dddot#1{\vbox{\ialign{##\crcr
      .\hskip-.5pt.\hskip-.5pt.\crcr\noalign{\kern1.5\p@\nointerlineskip}
      $\hfil\displaystyle{#1}\hfil$\crcr}}}

\newif\iftab@\tab@false
\newif\ifvtab@\vtab@false
\def\tab{\bgroup\tab@true\vtab@false\vst@bfalse\Strich@false%
   \def\\{\global\hline@@false%
     \ifhline@\global\hline@false\global\hline@@true\fi\cr}
   \edef\l@{\the\leftskip}\ialign\bgroup\hskip\l@##\hfil&&##\hfil\cr}
\def\endtab{\cr\egroup\egroup}
\def\vtab{\vtop\bgroup\vst@bfalse\vtab@true\tab@true\Strich@false%
   \bgroup\def\\{\cr}\ialign\bgroup&##\hfil\cr}
\def\endvtab{\cr\egroup\egroup\egroup}
\def\stab{\D@cke0.5pt\null 
 \bgroup\tab@true\vtab@false\vst@bfalse\Strich@true\Let@@\vspace@
 \normalbaselines\offinterlineskip
  \openup\spreadmlines@
 \edef\l@{\the\leftskip}\ialign
 \bgroup\hskip\l@##\hfil&&##\hfil\crcr}
\def\endstab{\crcr\egroup
 \egroup}
\newif\ifvst@b\vst@bfalse
\def\vstab{\D@cke0.5pt\null
 \vtop\bgroup\tab@true\vtab@false\vst@btrue\Strich@true\bgroup\Let@@\vspace@
 \normalbaselines\offinterlineskip
  \openup\spreadmlines@\bgroup}
\def\endvstab{\crcr\egroup\egroup
 \egroup\tab@false\Strich@false}

\newdimen\htstrut@
\htstrut@8.5\p@
\newdimen\htStrut@
\htStrut@12\p@
\newdimen\dpstrut@
\dpstrut@3.5\p@
\newdimen\dpStrut@
\dpStrut@3.5\p@
\def\openup{\afterassignment\@penup\dimen@=}
\def\@penup{\advance\lineskip\dimen@
  \advance\baselineskip\dimen@
  \advance\lineskiplimit\dimen@
  \divide\dimen@ by2
  \advance\htstrut@\dimen@
  \advance\htStrut@\dimen@
  \advance\dpstrut@\dimen@
  \advance\dpStrut@\dimen@}
\def\Let@@{\relax%
    \def\\{\global\hline@@false%
     \ifhline@\global\hline@false\global\hline@@true\fi\cr}%
    \iffalse}\fi}
\def\matrix{\null\,\vcenter\bgroup
 \tab@false\vtab@false\vst@bfalse\Strich@false\Let@@\vspace@
 \normalbaselines\openup\spreadmlines@\ialign
 \bgroup\hfil$\m@th##$\hfil&&\quad\hfil$\m@th##$\hfil\crcr
 \Mathstrut@\crcr\noalign{\kern-\baselineskip}}
\def\endmatrix{\crcr\Mathstrut@\crcr\noalign{\kern-\baselineskip}\egroup
 \egroup\,}
\def\smatrix{\D@cke0.5pt\null\,
 \vcenter\bgroup\tab@false\vtab@false\vst@bfalse\Strich@true\Let@@\vspace@
 \normalbaselines\offinterlineskip
  \openup\spreadmlines@\ialign
 \bgroup\hfil$\m@th##$\hfil&&\quad\hfil$\m@th##$\hfil\crcr}
\def\endsmatrix{\crcr\egroup
 \egroup\,\Strich@false}
\newdimen\D@cke
\def\Dicke#1{\global\D@cke#1}
\newtoks\tabs@\tabs@{&}
\newif\ifStrich@\Strich@false
\newif\iff@rst

\def\Stricherr@{\iftab@\ifvtab@\errmessage{\noexpand\s not allowed
     here. Use \noexpand\vstab!}%
  \else\errmessage{\noexpand\s not allowed here. Use \noexpand\stab!}%
  \fi\else\errmessage{\noexpand\s not allowed
     here. Use \noexpand\smatrix!}\fi}
\def\format{\ifvst@b\else\crcr\fi\egroup\iffalse{\fi\ifnum`}=0 \fi\format@}
\def\format@#1\\{\def\preamble@{#1}%
 \def\Str@chfehlt##1{\ifx##1\s\Stricherr@\fi\ifx##1\\\let\Next\relax%
   \else\let\Next\Str@chfehlt\fi\Next}%
 \def\c{\hfil\noexpand\ifhline@@\hbox{\vrule height\htStrut@%
   depth\dpstrut@ width\z@}\noexpand\fi%
   \ifStrich@\hbox{\vrule height\htstrut@ depth\dpstrut@ width\z@}%
   \fi\iftab@\else$\m@th\fi\the\hashtoks@\iftab@\else$\fi\hfil}%
 \def\r{\hfil\noexpand\ifhline@@\hbox{\vrule height\htStrut@%
   depth\dpstrut@ width\z@}\noexpand\fi%
   \ifStrich@\hbox{\vrule height\htstrut@ depth\dpstrut@ width\z@}%
   \fi\iftab@\else$\m@th\fi\the\hashtoks@\iftab@\else$\fi}%
 \def\l{\noexpand\ifhline@@\hbox{\vrule height\htStrut@%
   depth\dpstrut@ width\z@}\noexpand\fi%
   \ifStrich@\hbox{\vrule height\htstrut@ depth\dpstrut@ width\z@}%
   \fi\iftab@\else$\m@th\fi\the\hashtoks@\iftab@\else$\fi\hfil}%
 \def\s{\ifStrich@\ \the\tabs@\vrule width\D@cke\the\hashtoks@%
          \fi\the\tabs@\ }%
 \def\sa{\ifStrich@\vrule width\D@cke\the\hashtoks@%
            \the\tabs@\ %
            \fi}%
 \def\se{\ifStrich@\ \the\tabs@\vrule width\D@cke\the\hashtoks@\fi}%
 \def\cd{\hfil\noexpand\ifhline@@\hbox{\vrule height\htStrut@%
   depth\dpstrut@ width\z@}\noexpand\fi%
   \ifStrich@\hbox{\vrule height\htstrut@ depth\dpstrut@ width\z@}%
   \fi$\dsize\m@th\the\hashtoks@$\hfil}%
 \def\rd{\hfil\noexpand\ifhline@@\hbox{\vrule height\htStrut@%
   depth\dpstrut@ width\z@}\noexpand\fi%
   \ifStrich@\hbox{\vrule height\htstrut@ depth\dpstrut@ width\z@}%
   \fi$\dsize\m@th\the\hashtoks@$}%
 \def\ld{\noexpand\ifhline@@\hbox{\vrule height\htStrut@%
   depth\dpstrut@ width\z@}\noexpand\fi%
   \ifStrich@\hbox{\vrule height\htstrut@ depth\dpstrut@ width\z@}%
   \fi$\dsize\m@th\the\hashtoks@$\hfil}%
 \ifStrich@\else\Str@chfehlt#1\\\fi%
 \setbox\z@\hbox{\xdef\Preamble@{\preamble@}}\ifnum`{=0 \fi\iffalse}\fi
 \ialign\bgroup\span\Preamble@\crcr}
\newif\ifhline@\hline@false
\newif\ifhline@@\hline@@false
\def\hlinefor#1{\multispan@{\strip@#1 }\leaders\hrule height\D@cke\hfill%
    \global\hline@true\ignorespaces}
\def\Item "#1"{\par\noindent\hangindent2\parindent%
  \hangafter1\setbox0\hbox{\rm#1\enspace}\ifdim\wd0>2\parindent%
  \box0\else\hbox to 2\parindent{\rm#1\hfil}\fi\ignorespaces}
\def\ITEM #1"#2"{\par\noindent\hangafter1\hangindent#1%
  \setbox0\hbox{\rm#2\enspace}\ifdim\wd0>#1%
  \box0\else\hbox to 0pt{\rm#2\hss}\hskip#1\fi\ignorespaces}
\def\item"#1"{\par\noindent\hang%
  \setbox0=\hbox{\rm#1\enspace}\ifdim\wd0>\the\parindent%
  \box0\else\hbox to \parindent{\rm#1\hfil}\enspace\fi\ignorespaces}
\let\plainitem@\item
\catcode`\@=13

\catcode`\@=11
\font\tenln    = line10
\font\tenlnw   = linew10

\newskip\Einheit \Einheit=0.5cm
\newcount\xcoord \newcount\ycoord
\newdimen\xdim \newdimen\ydim \newdimen\PfadD@cke \newdimen\Pfadd@cke

\newcount\@tempcnta
\newcount\@tempcntb

\newdimen\@tempdima
\newdimen\@tempdimb

\newdimen\@wholewidth
\newdimen\@halfwidth

\newcount\@xarg
\newcount\@yarg
\newcount\@yyarg
\newbox\@linechar
\newbox\@tempboxa
\newdimen\@linelen
\newdimen\@clnwd
\newdimen\@clnht

\newif\if@negarg

\def\@whilenoop#1{}
\def\@whiledim#1\do #2{\ifdim #1\relax#2\@iwhiledim{#1\relax#2}\fi}
\def\@iwhiledim#1{\ifdim #1\let\@nextwhile=\@iwhiledim
        \else\let\@nextwhile=\@whilenoop\fi\@nextwhile{#1}}

\def\@whileswnoop#1\fi{}
\def\@whilesw#1\fi#2{#1#2\@iwhilesw{#1#2}\fi\fi}
\def\@iwhilesw#1\fi{#1\let\@nextwhile=\@iwhilesw
         \else\let\@nextwhile=\@whileswnoop\fi\@nextwhile{#1}\fi}

\def\thinlines{\let\@linefnt\tenln \let\@circlefnt\tencirc
  \@wholewidth\fontdimen8\tenln \@halfwidth .5\@wholewidth}
\def\thicklines{\let\@linefnt\tenlnw \let\@circlefnt\tencircw
  \@wholewidth\fontdimen8\tenlnw \@halfwidth .5\@wholewidth}
\thinlines

\PfadD@cke1pt \Pfadd@cke0.5pt
\def\PfadDicke#1{\PfadD@cke#1 \divide\PfadD@cke by2 \Pfadd@cke\PfadD@cke \multiply\PfadD@cke by2}
\long\def\LOOP#1\REPEAT{\def\BODY{#1}\ITERATE}
\def\ITERATE{\BODY \let\next\ITERATE \else\let\next\relax\fi \next}
\let\REPEAT=\fi
\def\Punkt{\hbox{\raise-2pt\hbox to0pt{\hss$\ssize\bullet$\hss}}}
\def\DuennPunkt(#1,#2){\unskip
  \raise#2 \Einheit\hbox to0pt{\hskip#1 \Einheit
          \raise-2.5pt\hbox to0pt{\hss$\bullet$\hss}\hss}}
\def\NormalPunkt(#1,#2){\unskip
  \raise#2 \Einheit\hbox to0pt{\hskip#1 \Einheit
          \raise-3pt\hbox to0pt{\hss\twelvepoint$\bullet$\hss}\hss}}
\def\DickPunkt(#1,#2){\unskip
  \raise#2 \Einheit\hbox to0pt{\hskip#1 \Einheit
          \raise-4pt\hbox to0pt{\hss\fourteenpoint$\bullet$\hss}\hss}}
\def\Kreis(#1,#2){\unskip
  \raise#2 \Einheit\hbox to0pt{\hskip#1 \Einheit
          \raise-4pt\hbox to0pt{\hss\fourteenpoint$\circ$\hss}\hss}}

\def\Line@(#1,#2)#3{\@xarg #1\relax \@yarg #2\relax
\@linelen=#3\Einheit
\ifnum\@xarg =0 \@vline
  \else \ifnum\@yarg =0 \@hline \else \@sline\fi
\fi}

\def\@sline{\ifnum\@xarg< 0 \@negargtrue \@xarg -\@xarg \@yyarg -\@yarg
  \else \@negargfalse \@yyarg \@yarg \fi
\ifnum \@yyarg >0 \@tempcnta\@yyarg \else \@tempcnta -\@yyarg \fi
\ifnum\@tempcnta>6 \@badlinearg\@tempcnta0 \fi
\ifnum\@xarg>6 \@badlinearg\@xarg 1 \fi
\setbox\@linechar\hbox{\@linefnt\@getlinechar(\@xarg,\@yyarg)}%
\ifnum \@yarg >0 \let\@upordown\raise \@clnht\z@
   \else\let\@upordown\lower \@clnht \ht\@linechar\fi
\@clnwd=\wd\@linechar
\if@negarg \hskip -\wd\@linechar \def\@tempa{\hskip -2\wd\@linechar}\else
     \let\@tempa\relax \fi
\@whiledim \@clnwd <\@linelen \do
  {\@upordown\@clnht\copy\@linechar
   \@tempa
   \advance\@clnht \ht\@linechar
   \advance\@clnwd \wd\@linechar}%
\advance\@clnht -\ht\@linechar
\advance\@clnwd -\wd\@linechar
\@tempdima\@linelen\advance\@tempdima -\@clnwd
\@tempdimb\@tempdima\advance\@tempdimb -\wd\@linechar
\if@negarg \hskip -\@tempdimb \else \hskip \@tempdimb \fi
\multiply\@tempdima \@m
\@tempcnta \@tempdima \@tempdima \wd\@linechar \divide\@tempcnta \@tempdima
\@tempdima \ht\@linechar \multiply\@tempdima \@tempcnta
\divide\@tempdima \@m
\advance\@clnht \@tempdima
\ifdim \@linelen <\wd\@linechar
   \hskip \wd\@linechar
  \else\@upordown\@clnht\copy\@linechar\fi}

\def\@hline{\ifnum \@xarg <0 \hskip -\@linelen \fi
\vrule height\Pfadd@cke width \@linelen depth\Pfadd@cke
\ifnum \@xarg <0 \hskip -\@linelen \fi}

\def\@getlinechar(#1,#2){\@tempcnta#1\relax\multiply\@tempcnta 8
\advance\@tempcnta -9 \ifnum #2>0 \advance\@tempcnta #2\relax\else
\advance\@tempcnta -#2\relax\advance\@tempcnta 64 \fi
\char\@tempcnta}

\def\Vektor(#1,#2)#3(#4,#5){\unskip\leavevmode
  \xcoord#4\relax \ycoord#5\relax
      \raise\ycoord \Einheit\hbox to0pt{\hskip\xcoord \Einheit
         \Vector@(#1,#2){#3}\hss}}

\def\Vector@(#1,#2)#3{\@xarg #1\relax \@yarg #2\relax
\@tempcnta \ifnum\@xarg<0 -\@xarg\else\@xarg\fi
\ifnum\@tempcnta<5\relax
\@linelen=#3\Einheit
\ifnum\@xarg =0 \@vvector
  \else \ifnum\@yarg =0 \@hvector \else \@svector\fi
\fi
\else\@badlinearg\fi}

\def\@hvector{\@hline\hbox to 0pt{\@linefnt
\ifnum \@xarg <0 \@getlarrow(1,0)\hss\else
    \hss\@getrarrow(1,0)\fi}}

\def\@vvector{\ifnum \@yarg <0 \@downvector \else \@upvector \fi}

\def\@svector{\@sline
\@tempcnta\@yarg \ifnum\@tempcnta <0 \@tempcnta=-\@tempcnta\fi
\ifnum\@tempcnta <5
  \hskip -\wd\@linechar
  \@upordown\@clnht \hbox{\@linefnt  \if@negarg
  \@getlarrow(\@xarg,\@yyarg) \else \@getrarrow(\@xarg,\@yyarg) \fi}%
\else\@badlinearg\fi}

\def\@upline{\hbox to \z@{\hskip -.5\Pfadd@cke \vrule width \Pfadd@cke
   height \@linelen depth \z@\hss}}

\def\@downline{\hbox to \z@{\hskip -.5\Pfadd@cke \vrule width \Pfadd@cke
   height \z@ depth \@linelen \hss}}

\def\@upvector{\@upline\setbox\@tempboxa\hbox{\@linefnt\char'66}\raise
     \@linelen \hbox to\z@{\lower \ht\@tempboxa\box\@tempboxa\hss}}

\def\@downvector{\@downline\lower \@linelen
      \hbox to \z@{\@linefnt\char'77\hss}}

\def\@getlarrow(#1,#2){\ifnum #2 =\z@ \@tempcnta='33\else
\@tempcnta=#1\relax\multiply\@tempcnta \sixt@@n \advance\@tempcnta
-9 \@tempcntb=#2\relax\multiply\@tempcntb \tw@
\ifnum \@tempcntb >0 \advance\@tempcnta \@tempcntb\relax
\else\advance\@tempcnta -\@tempcntb\advance\@tempcnta 64
\fi\fi\char\@tempcnta}

\def\@getrarrow(#1,#2){\@tempcntb=#2\relax
\ifnum\@tempcntb < 0 \@tempcntb=-\@tempcntb\relax\fi
\ifcase \@tempcntb\relax \@tempcnta='55 \or
\ifnum #1<3 \@tempcnta=#1\relax\multiply\@tempcnta
24 \advance\@tempcnta -6 \else \ifnum #1=3 \@tempcnta=49
\else\@tempcnta=58 \fi\fi\or
\ifnum #1<3 \@tempcnta=#1\relax\multiply\@tempcnta
24 \advance\@tempcnta -3 \else \@tempcnta=51\fi\or
\@tempcnta=#1\relax\multiply\@tempcnta
\sixt@@n \advance\@tempcnta -\tw@ \else
\@tempcnta=#1\relax\multiply\@tempcnta
\sixt@@n \advance\@tempcnta 7 \fi\ifnum #2<0 \advance\@tempcnta 64 \fi
\char\@tempcnta}

\def\Diagonale(#1,#2)#3{\unskip\leavevmode
  \xcoord#1\relax \ycoord#2\relax
      \raise\ycoord \Einheit\hbox to0pt{\hskip\xcoord \Einheit
         \Line@(1,1){#3}\hss}}
\def\AntiDiagonale(#1,#2)#3{\unskip\leavevmode
  \xcoord#1\relax \ycoord#2\relax 
      \raise\ycoord \Einheit\hbox to0pt{\hskip\xcoord \Einheit
         \Line@(1,-1){#3}\hss}}
\def\Pfad(#1,#2),#3\endPfad{\unskip\leavevmode
  \xcoord#1 \ycoord#2 \thicklines\ZeichnePfad#3\endPfad\thinlines}
\def\ZeichnePfad#1{\ifx#1\endPfad\let\next\relax
  \else\let\next\ZeichnePfad
    \ifnum#1=1
      \raise\ycoord \Einheit\hbox to0pt{\hskip\xcoord \Einheit
         \vrule height\Pfadd@cke width1 \Einheit depth\Pfadd@cke\hss}%
      \advance\xcoord by 1
    \else\ifnum#1=2
      \raise\ycoord \Einheit\hbox to0pt{\hskip\xcoord \Einheit
        \hbox{\hskip-\PfadD@cke\vrule height1 \Einheit width\PfadD@cke depth0pt}\hss}%
      \advance\ycoord by 1
    \else\ifnum#1=3
      \raise\ycoord \Einheit\hbox to0pt{\hskip\xcoord \Einheit
         \Line@(1,1){1}\hss}
      \advance\xcoord by 1
      \advance\ycoord by 1
    \else\ifnum#1=4
      \raise\ycoord \Einheit\hbox to0pt{\hskip\xcoord \Einheit
         \Line@(1,-1){1}\hss}
      \advance\xcoord by 1
      \advance\ycoord by -1
    \else\ifnum#1=5
      \advance\xcoord by -1
      \raise\ycoord \Einheit\hbox to0pt{\hskip\xcoord \Einheit
         \vrule height\Pfadd@cke width1 \Einheit depth\Pfadd@cke\hss}%
    \else\ifnum#1=6
      \advance\ycoord by -1
      \raise\ycoord \Einheit\hbox to0pt{\hskip\xcoord \Einheit
        \hbox{\hskip-\PfadD@cke\vrule height1 \Einheit width\PfadD@cke depth0pt}\hss}%
    \else\ifnum#1=7
      \advance\xcoord by -1
      \advance\ycoord by -1
      \raise\ycoord \Einheit\hbox to0pt{\hskip\xcoord \Einheit
         \Line@(1,1){1}\hss}
    \else\ifnum#1=8
      \advance\xcoord by -1
      \advance\ycoord by +1
      \raise\ycoord \Einheit\hbox to0pt{\hskip\xcoord \Einheit
         \Line@(1,-1){1}\hss}
    \fi\fi\fi\fi
    \fi\fi\fi\fi
  \fi\next}
\def\hSSchritt{\leavevmode\raise-.4pt\hbox to0pt{\hss.\hss}\hskip.2\Einheit
  \raise-.4pt\hbox to0pt{\hss.\hss}\hskip.2\Einheit
  \raise-.4pt\hbox to0pt{\hss.\hss}\hskip.2\Einheit
  \raise-.4pt\hbox to0pt{\hss.\hss}\hskip.2\Einheit
  \raise-.4pt\hbox to0pt{\hss.\hss}\hskip.2\Einheit}
\def\vSSchritt{\vbox{\baselineskip.2\Einheit\lineskiplimit0pt
\hbox{.}\hbox{.}\hbox{.}\hbox{.}\hbox{.}}}
\def\DSSchritt{\leavevmode\raise-.4pt\hbox to0pt{%
  \hbox to0pt{\hss.\hss}\hskip.2\Einheit
  \raise.2\Einheit\hbox to0pt{\hss.\hss}\hskip.2\Einheit
  \raise.4\Einheit\hbox to0pt{\hss.\hss}\hskip.2\Einheit
  \raise.6\Einheit\hbox to0pt{\hss.\hss}\hskip.2\Einheit
  \raise.8\Einheit\hbox to0pt{\hss.\hss}\hss}}
\def\dSSchritt{\leavevmode\raise-.4pt\hbox to0pt{%
  \hbox to0pt{\hss.\hss}\hskip.2\Einheit
  \raise-.2\Einheit\hbox to0pt{\hss.\hss}\hskip.2\Einheit
  \raise-.4\Einheit\hbox to0pt{\hss.\hss}\hskip.2\Einheit
  \raise-.6\Einheit\hbox to0pt{\hss.\hss}\hskip.2\Einheit
  \raise-.8\Einheit\hbox to0pt{\hss.\hss}\hss}}
\def\SPfad(#1,#2),#3\endSPfad{\unskip\leavevmode
  \xcoord#1 \ycoord#2 \ZeichneSPfad#3\endSPfad}
\def\ZeichneSPfad#1{\ifx#1\endSPfad\let\next\relax
  \else\let\next\ZeichneSPfad
    \ifnum#1=1
      \raise\ycoord \Einheit\hbox to0pt{\hskip\xcoord \Einheit
         \hSSchritt\hss}%
      \advance\xcoord by 1
    \else\ifnum#1=2
      \raise\ycoord \Einheit\hbox to0pt{\hskip\xcoord \Einheit
        \hbox{\hskip-2pt \vSSchritt}\hss}%
      \advance\ycoord by 1
    \else\ifnum#1=3
      \raise\ycoord \Einheit\hbox to0pt{\hskip\xcoord \Einheit
         \DSSchritt\hss}
      \advance\xcoord by 1
      \advance\ycoord by 1
    \else\ifnum#1=4
      \raise\ycoord \Einheit\hbox to0pt{\hskip\xcoord \Einheit
         \dSSchritt\hss}
      \advance\xcoord by 1
      \advance\ycoord by -1
    \else\ifnum#1=5
      \advance\xcoord by -1
      \raise\ycoord \Einheit\hbox to0pt{\hskip\xcoord \Einheit
         \hSSchritt\hss}%
    \else\ifnum#1=6
      \advance\ycoord by -1
      \raise\ycoord \Einheit\hbox to0pt{\hskip\xcoord \Einheit
        \hbox{\hskip-2pt \vSSchritt}\hss}%
    \else\ifnum#1=7
      \advance\xcoord by -1
      \advance\ycoord by -1
      \raise\ycoord \Einheit\hbox to0pt{\hskip\xcoord \Einheit
         \DSSchritt\hss}
    \else\ifnum#1=8
      \advance\xcoord by -1
      \advance\ycoord by 1
      \raise\ycoord \Einheit\hbox to0pt{\hskip\xcoord \Einheit
         \dSSchritt\hss}
    \fi\fi\fi\fi
    \fi\fi\fi\fi
  \fi\next}
\def\Koordinatenachsen(#1,#2){\unskip
 \hbox to0pt{\hskip-.5pt\vrule height#2 \Einheit width.5pt depth1 \Einheit}%
 \hbox to0pt{\hskip-1 \Einheit \xcoord#1 \advance\xcoord by1
    \vrule height0.25pt width\xcoord \Einheit depth0.25pt\hss}}
\def\Koordinatenachsen(#1,#2)(#3,#4){\unskip
 \hbox to0pt{\hskip-.5pt \ycoord-#4 \advance\ycoord by1
    \vrule height#2 \Einheit width.5pt depth\ycoord \Einheit}%
 \hbox to0pt{\hskip-1 \Einheit \hskip#3\Einheit 
    \xcoord#1 \advance\xcoord by1 \advance\xcoord by-#3 
    \vrule height0.25pt width\xcoord \Einheit depth0.25pt\hss}}
\def\Gitter(#1,#2){\unskip \xcoord0 \ycoord0 \leavevmode
  \LOOP\ifnum\ycoord<#2
    \loop\ifnum\xcoord<#1
      \raise\ycoord \Einheit\hbox to0pt{\hskip\xcoord \Einheit\Punkt\hss}%
      \advance\xcoord by1
    \repeat
    \xcoord0
    \advance\ycoord by1
  \REPEAT}
\def\Gitter(#1,#2)(#3,#4){\unskip \xcoord#3 \ycoord#4 \leavevmode
  \LOOP\ifnum\ycoord<#2
    \loop\ifnum\xcoord<#1
      \raise\ycoord \Einheit\hbox to0pt{\hskip\xcoord \Einheit\Punkt\hss}%
      \advance\xcoord by1
    \repeat
    \xcoord#3
    \advance\ycoord by1
  \REPEAT}
\def\Label#1#2(#3,#4){\unskip \xdim#3 \Einheit \ydim#4 \Einheit
  \def\lo{\advance\xdim by-.5 \Einheit \advance\ydim by.5 \Einheit}%
  \def\llo{\advance\xdim by-.25cm \advance\ydim by.5 \Einheit}%
  \def\loo{\advance\xdim by-.5 \Einheit \advance\ydim by.25cm}%
  \def\o{\advance\ydim by.25cm}%
  \def\ro{\advance\xdim by.5 \Einheit \advance\ydim by.5 \Einheit}%
  \def\rro{\advance\xdim by.25cm \advance\ydim by.5 \Einheit}%
  \def\roo{\advance\xdim by.5 \Einheit \advance\ydim by.25cm}%
  \def\l{\advance\xdim by-.30cm}%
  \def\r{\advance\xdim by.30cm}%
  \def\lu{\advance\xdim by-.5 \Einheit \advance\ydim by-.6 \Einheit}%
  \def\llu{\advance\xdim by-.25cm \advance\ydim by-.6 \Einheit}%
  \def\luu{\advance\xdim by-.5 \Einheit \advance\ydim by-.30cm}%
  \def\u{\advance\ydim by-.30cm}%
  \def\ru{\advance\xdim by.5 \Einheit \advance\ydim by-.6 \Einheit}%
  \def\rru{\advance\xdim by.25cm \advance\ydim by-.6 \Einheit}%
  \def\ruu{\advance\xdim by.5 \Einheit \advance\ydim by-.30cm}%
  #1\raise\ydim\hbox to0pt{\hskip\xdim
     \vbox to0pt{\vss\hbox to0pt{\hss$#2$\hss}\vss}\hss}%
}
\catcode`\@=13

\hsize13cm
\vsize19cm
\newdimen\fullhsize
\newdimen\fullvsize
\newdimen\halfsize
\fullhsize13cm
\fullvsize19cm
\halfsize=0.5\fullhsize
\advance\halfsize by-0.5em

\magnification1200

\TagsOnRight

\def\AignAA{1}
\def\AignAD{2}
\def\AignAB{3}
\def\BeCaAA{4}
\def\BeCQAA{5}
\def\BeCQAB{6}
\def\BrEsAB{7}
\def\BrEsAC{8}
\def\BresAO{9}
\def\ChFrAA{10}
\def\ChGoAA{11}
\def\CiglAM{12}
\def\CiglAS{13}
\def\CiKrAA{14}
\def\CvRIAA{15}
\def\DeViAB{16}
\def\DesaAB{17}
\def\EgRRAA{18}
\def\EgRRAB{19}
\def\FishAA{20}
\def\GeViAA{21}
\def\GeViAB{22}
\def\GeXiAA{23}
\def\GhKrAA{24}
\def\GrKPAA{25}
\def\GrJSAA{26}
\def\HiPeAA{27}
\def\JoSaAB{28}
\def\KratAM{29}
\def\KratBN{30}
\def\KratBZ{31}
\def\KratBW{32}
\def\LaymAA{33}
\def\LindAA{34}
\def\MaWoAA{35}
\def\MohaAE{36}
\def\OwEBAA{37}
\def\PearAA{38}
\def\RadoAI{39}
\def\RadoAB{40}
\def\RadoAH{41}
\def\RadoAG{42}
\def\StanBI{43}
\def\TammAD{44}
\def\VienAE{45}
\def\VienAB{46}

\def\AA{1.5}
\def\AB{1.1}
\def\AC{1.2}
\def\AD{1.3}
\def\AE{1.4}
\def\AF{1.6}
\def\AG{1.7}
\def\AH{3.1}
\def\AI{3.2}
\def\AJ{4.1}
\def\AJb{4.2}
\def\AJa{4.3}
\def\AK{4.4}
\def\AL{4.5}
\def\ALa{4.6}
\def\AM{5.1}
\def\AN{5.2}
\def\AO{5.3}
\def\AP{5.4}
\def\AQ{5.5}
\def\AR{5.6}
\def\AS{5.7}
\def\AT{5.8}
\def\AU{5.9}
\def\AV{5.10}
\def\AW{5.11}

\def\TA{1}
\def\TAA{2}
\def\TB{3}
\def\TC{4}
\def\TD{5}
\def\TE{6}
\def\TF{7}
\def\TG{8}
\def\TH{9}
\def\TI{10}

\def\FA{1}
\def\FB{2}
\def\FC{3}
\def\FD{4}
\def\FE{5}

\topmatter 
\title Determinants of (generalised) Catalan numbers
\endtitle 
\author C.~Krattenthaler$^{\dagger}$
\endauthor 
\affil 
Fakult\"at f\"ur Mathematik, Universit\"at Wien,\\
Nordbergstra{\ss}e~15, A-1090 Vienna, Austria.\\
WWW: \tt http://www.mat.univie.ac.at/\~{}kratt
\endaffil
\address Fakult\"at f\"ur Mathematik, Universit\"at Wien,
Nordbergstra{\ss}e~15, A-1090 Vienna, Austria.
WWW: \tt http://www.mat.univie.ac.at/\~{}kratt
\endaddress
\thanks $^\dagger$Research partially supported 
by the Austrian Science Foundation FWF, grants Z130-N13 and S9607-N13,
the latter in the framework of the National Research Network
``Analytic Combinatorics and Probabilistic Number Theory"%
\endthanks
\subjclass Primary 05A19;
 Secondary 05A10 05A15 05E10 11C20 15A15 33C45 
\endsubjclass
\keywords Hankel determinants, generalised Catalan numbers, Fibonacci numbers,
Dodg\-son condensation, non-intersecting lattice paths
\endkeywords
\abstract 
We show that recent determinant evaluations involving Catalan numbers
and generalisations thereof have most convenient explanations by
combining the Lindstr\"om--Gessel--Viennot theorem on non-intersecting 
lattice paths with a simple determinant lemma from
[{\it Manuscripta Math\.} {\bf 69} (1990), 173--202].
This approach leads also naturally to extensions and
generalisations.
\endabstract
\endtopmatter
\document

\subhead 1. Introduction\endsubhead
Determinants of matrices containing combinatorial numbers have always
been of big attraction to many researchers. The combinatorial numbers
which are in the centre of the present paper are the
{\it Catalan numbers} $C_n$, defined by 
$$C_n=\frac {1} {n+1}\binom {2n}n,$$
and extensions thereof. (The reader is referred to \cite{\StanBI, 
Exercise~6.19} for extensive information on these numbers.)
Numerous papers exist in the literature featuring determinants of matrices
the entries of which contain Catalan numbers or generalisations thereof, see
\cite{\AignAA, \AignAD, \AignAB, \BeCaAA, \BeCQAA, \BeCQAB, \BrEsAB, \BrEsAC,
\ChFrAA, \ChGoAA, \CiglAM, \CiglAS, 
\CiKrAA, \CvRIAA, \DeViAB, \EgRRAA, \GeViAB, \GeXiAA, \GhKrAA, \KratBW, \LaymAA,
\MaWoAA, \OwEBAA, \PearAA, \RadoAI, \RadoAB, \RadoAH, \RadoAG, 
\TammAD, \VienAE} (and this list is for sure incomplete). 
The contexts in which they appear are also
widespread, ranging from lattice path enumeration (cf\.
\cite{\AignAD, \AignAB, \BeCQAA, \BeCQAB, \BrEsAB, \BrEsAC, 
\CiglAM, \CiglAS, \GeXiAA, \EgRRAA,
\KratBW, \MaWoAA, \OwEBAA, \VienAE}), 
plane partition
and tableaux counting (cf\. \cite{\AignAD, \ChGoAA, \DeViAB, \GeViAB}), 
continued fractions and orthogonal polynomials 
(cf\. \cite{\CiglAS, \CvRIAA, \RadoAB, \TammAD, \VienAE}), 
statistical physics (cf\. \cite{\BrEsAB, \BrEsAC, \KratBW, \OwEBAA}),
up to commutative
algebra and algebraic geometry (cf\. \cite{\GhKrAA}).

One of the most popular themes in this context is {\it Hankel
determinants} of Catalan numbers, that is, determinants of the form 
$\det_{0\le i,j\le n-1}(a_{i+j})$, where the sequence $(a_i)_{i\ge0}$
involves Catalan numbers.  (Cf\. \cite{\KratBN, Sec.~2.7} and
\cite{\KratBZ, Sec.~5.4} for general information on the evaluation of
Hankel determinants.) To be more specific, Hankel determinant evaluations
such as 
$$\align 
\det_{0\le i,j\le n-1}(C_{i+j})&=1,
\tag\AB\\
\det_{0\le i,j\le n-1}(C_{i+j+1})&=1,
\tag\AC
\endalign$$
or
$$\det_{0\le i,j\le n-1}(C_{i+j+2})=n+1
\tag\AD$$
have been addressed numerous times in the literature.
More recently, it has been observed in \cite{\LaymAA} and proved in
\cite{\CvRIAA} that
$$\det_{0\le i,j\le n-1}(C_{i+j}+C_{i+j+1})=F_{2n},
\tag\AE$$
where $F_m$ denotes the $m$-th {\it Fibonacci number}; that is,
$F_0=F_1=1$ and $F_m=F_{m-1}+F_{m-2}$ for $m\ge2$.

Let $k$ be a fixed positive integer.
Taking the {\it generalised Catalan number} $C_{n,k}$ from \cite{\HiPeAA},
given by
$$C_{n,k}=\frac {n-(k-1)\lfloor{\frac {n} {k-1}}\rfloor+1} 
{n+\lfloor{\frac {n}
{k-1}}\rfloor+1}
\binom {n+\lfloor{\frac {n} {k-1}}\rfloor+1} {n+1},
\tag\AA
$$
the determinant evaluation in (\AE) has been generalised in \cite{\ChFrAA}
to
$$\det_{0\le i,j\le n-1}(C_{(k-1)i+j,k}+C_{(k-1)(i+1)+j,k})=
\sum _{s=0} ^{n}\binom {(k-1)s+n} {n-s}
\tag\AF$$
and
$$\det_{0\le i,j\le n-1}(C_{(k-1)i+j,k}+C_{(k-1)i+j+1,k})=
\sum _{s=0} ^{n}\binom {\lfloor\frac {s} {k-1}\rfloor+n} {n-s}.
\tag\AG$$

Many different methods have been employed to prove the formulae
(\AB)--(\AE), (\AF), and (\AG), among which are
direct determinant manipulations,
non-intersecting lattice paths, orthogonal polynomials, and LU
factorisation. However, in the author's opinion, none of the
published approaches explains in a satisfactory and uniform way 
why these identities exist. This is what we aim to do in this
article. We show that the above determinant evaluations are actually
special cases of families of more general determinant evaluations,
which result from a combination of a simple determinant lemma from
\cite{\KratAM} (see Lemma~\TA) and the main theorem on
non-intersecting lattice paths (see Theorem~\TAA). In particular,
Theorem~\TB, originally due to Gessel and Viennot \cite{\GeViAB}, 
shows that one can in fact evaluate the determinant
$\det_{0\le i,j\le n-1}(C_{\al_i+j})$ in closed form, where
$\al_0,\al_1,\dots,\al_{n-1}$ are arbitrary non-negative numbers,
thus largely generalising (\AB)--(\AD).
Together with the simple determinant expansion in Lemma~\TC, this
explains as well formula~(\AE). Moreover, it also enables us to embed (\AE)
in a much larger family of determinant evaluations, see
Corollary~\TD. 

Starting point for our ``explanation" of (\AF) and (\AG) is the
(well-known) observation that the generalised
Catalan numbers in (\AA) count lattice paths which stay below
a slanted line (see (\AJ)). Then, altogether, the main theorem on
non-intersecting lattice paths, combinatorial arguments involving
paths, and the aforementioned determinant lemma, lead in fact again to
much more general determinant evaluations, which we present in
Theorem~\TE\ (which generalises Theorem~\TB, and, thus,
identities~(\AB)--(\AD)), Corollary~\TF, and Theorem~\TH, respectively.

The final section, Section~5, is provided here for the sake of
completeness. There, we consider generalised Catalan
numbers which are different from those in (\AA). We briefly survey
the deep results from \cite{\EgRRAA, \GeXiAA}
which are known on Hankel determinants involving these numbers. 

We emphasise that we do not address {\it weighted\/} generalisations in this
paper. Sometimes these can also be approached by the method of
non-intersecting lattice paths (cf\. \cite{\CiglAS, \KratBW,
\OwEBAA, \VienAE} for examples), but often more refined methods are
needed (cf\. \cite{\BrEsAB, \BrEsAC, \CiglAM, \CiKrAA, \KratBW, \OwEBAA}). 

\subhead 2. Preliminaries\endsubhead
In this section we present the two auxiliary results on which
everything else in this paper is based: a general determinant lemma
and the Lindstr\"om--Gessel--Viennot theorem on non-intersecting
lattice paths.

We begin with the determinant lemma, which is taken from
\cite{\KratAM, Lemma~2.2} (see also \cite{\KratBN, Lemma~3}). 
There are many ways to prove it.
The ``easiest" is by condensation (once one takes the 
determinant identity by Desnanot and Jacobi which underlies the
condensation method for granted; cf\. \cite{\KratBN,
Sec.~2.3}). For alternative proofs see \cite{\BresAO, Theorem~2.9},
\cite{\KratAM}, and \cite{\KratBN, App.~B}.

\proclaim{Lemma~\TA}
Let $X_0,\dots,X_{n-1}$, $A_1,\dots,A_{n-1}$, and $B_1,
\dots ,B_{n-1}$ be indeterminates.\linebreak 
Then there holds
$$\multline 
\det_{0\le i,j\le n-1}\Big((X_i+A_{n-1})(X_i+A_{n-2})\cdots(X_i+A_{j+1})
(X_i+B_{j})(X_i+B_{j-1})\cdots (X_i+B_1)\Big)\\
\hskip2cm =\prod _{0\le i<j\le n-1} ^{}(X_i-X_j)
\prod _{1\le i\le j\le n-1} ^{}(B_i-A_j). 
\endmultline$$
\endproclaim

The other auxiliary result that we need is a determinant formula for
the enumeration of non-intersecting lattice paths, originally due to
Lindstr\"om \cite{\LindAA, Lemma~1}, rediscovered by Gessel and Viennot
\cite{\GeViAB} (and special cases in \cite{\FishAA, Sec.~5.3} and in
\cite{\JoSaAB, \GrJSAA}; for a more detailed historical account
see Footnote~6 in \cite{\KratBW}). Since, as is the case in most
applications, we do not need the theorem in its most general form,
we state the special case that serves our purposes. The lattice paths
to which we want to apply the theorem are lattice paths in the
integer plane consisting of horizontal and vertical unit steps in the
positive direction. We make the convention for the rest of the paper
that this is what we mean whenever we speak of ``lattice paths,"
unless it is specified otherwise.
A family $(P_0,P_1,\dots,P_{n-1})$ of lattice paths $P_i$,
$i=0,1,\dots,n-1$, is called {\it non-intersecting} if
no two paths in the family have a point in common.
Examples can be found in Figure~\FA--\FE.
With the above terminology, we have the following theorem.

\proclaim{Theorem~\TAA}
Let $A_0,A_1,\dots,A_{n-1}$ and $E_0,E_1,\dots,E_{n-1}$ be lattice points
in the plane integer lattice
such that for $i<j$ and $k<l$ any lattice path from $A_i$ to $E_l$ has a
common point with any lattice path from $A_j$ to $E_k$.
Then the number of all families $(P_0,P_1,\dots,P_{n-1})$ of
non-intersecting lattice paths, $P_i$ running from $A_i$ to $E_i$,
$i=0,1,\dots,n-1$, is given by
$$\det_{0\le i,j\le n-1}\big(P(A_j\to E_i)\big),$$
where $P(A\to E)$ denotes the number of all lattice paths from $A$ to
$E$.
\endproclaim

\subhead 3. Determinants of Catalan numbers\endsubhead
In this section we present evaluations of 
determinants in which the entries involve
Catalan numbers. These determinant evaluations generalise Eqs.~(\AB)--(\AE). 
Our first theorem is a common generalisation of Eqs.~(\AB)--(\AD).
It is originally due to Gessel and Viennot \cite{\GeViAB, 
paragraph between Theorems~21 and 22}, who gave essentially 
the proof that we present below.

\proclaim{Theorem~\TB}
Let $n$ be a positive integer and $\al_0,\al_1,\dots,\al_{n-1}$
non-negative integers. Then
$$\det_{0\le i,j\le n-1}\big(C_{\al_i+j}\big)=
\prod _{0\le i<j\le n-1} ^{}(\al_j-\al_i)
\prod _{i=0} ^{n-1}\frac {(i+n)!\,(2\al_i)!} {(2i)!\,\al_i!\,(\al_i+n)!}.$$
\endproclaim
\demo{Proof}
We use the observation $C_m=(-1)^m2^{2m+1}\binom {1/2}{m+1}$ to
rewrite the determinant in question as
$$
(-1)^{\binom n2+\sum _{i=0} ^{n-1}\al_i}2^{n^2+2\sum _{i=0}
^{n-1}\al_i}\det_{0\le i,j\le n-1}\(\binom {1/2}{\al_i+j+1}\).
\tag\AH$$
Then, by taking some factors out of the $i$-th row of the matrix of
which we want to compute the determinant, $i=0,1,\dots,n-1$, we obtain
$$\multline 
\det_{0\le i,j\le n-1}\big(C_{\al_i+j}\big)=
(-1)^{\sum _{i=0} ^{n-1}\al_i}2^{n^2+2\sum _{i=0}
^{n-1}\al_i}\prod _{i=0} ^{n-1}\frac 
{\big(\frac {1} {2}-\al_i\big)_{\al_i+1}} {(\al_i+n)!}\\
\times
\det_{0\le i,j\le n-1}\bigg(\(\al_i+j-\frac {1} {2}\)
\(\al_i+j-\frac {3} {2}\)\cdots\(\al_i+\frac {1} {2}\)\\
\cdot
(\al_i+j+2)(\al_i+j+3)\cdots(\al_i+n)\bigg).
\endmultline$$
Now Lemma~\TA\ can be applied with $X_i=\al_i$, $A_j=j+1$, and
$B_j=j-\frac {1} {2}$. After some manipulation, one arrives at the
claimed result.\quad \quad \qed
\enddemo

\remark{Remarks} 
(a) Clearly, the determinant evaluations (\AB)--(\AD) are all special cases 
of Theorem~\TB.

\smallskip
(b) From a conceptual point of view, the ``actual" theorem is the determinant
evaluation for the determinant in (\AH). In fact, there holds the
determinant evaluation (see \cite{\KratBN, Theorem~26, (3.12)})
$$
\det_{0\le i,j\le n-1}\(\binom
A {\al_i+j}\)
=\frac {\prod _{0\le i<j\le n-1}
^{}(\al_i-\al_j)} {\prod _{i=0} ^{n-1}(\al_i+n-1)!}
\frac {\prod _{i=0} ^{n-1}(A+i)!} {\prod _{i=0} ^{n-1}(A-\al_i)!}.
$$
(The cited theorem in \cite{\KratBN} is even a
$q$-analogue of the above identity.)
This identity is not restricted to integral $A$ and $\al_i$: if one
interprets the appearing factorials and binomials as suitable 
expressions involving gamma functions
(cf.\ \cite{\GrKPAA, \S5.5, (5.96), (5.100)}), then
it continues to hold for real or complex $A$ and $\al_i$ as long as there do
not appear any singularities.

\smallskip
(c) Gessel and Viennot used Theorem~\TB\ to derive tableaux
enumeration results, see \cite{\GeViAB, Theorems~22 and 23}.

\smallskip
(d) The most elegant proof of the special cases (\AB)--(\AD) 
of Theorem~\TB\ is by
using non-intersecting lattice paths. For (\AB) and (\AC), these
(one figure) proofs appear as side results already in \cite{\VienAE}.
For (\AD), a corresponding proof is explained in
\cite{\BeCQAA}, as well as
a non-intersecting lattice path proof for the ``next" case, 
the evaluation of the Hankel determinant 
$\det_{0\le i,j\le n-1}\big(C_{i+j+3}\big)$. 
On the other hand, it seems unlikely that the full statement of
Theorem~\TB\ can be explained combinatorially. 

\smallskip
(e) An ubiquitous determinant is the Hankel determinant
$\det_{0\le i,j\le n-1}\big(C_{i+j+\be}\big)$, where $\be$ is some
fixed positive integer, which, by Theorem~\TB\ with $\al_i=i+\be$,
$i=0,1,\dots,n-1$, has a closed form product evaluation. 
By Theorem~\TAA\ and the standard
combinatorial interpretation of the Catalan number $C_m$ as the
number of lattice paths from the origin to $(m,m)$ which never pass
above the diagonal $x=y$, this determinant counts families
$(P_0,P_1,\dots,P_{n-1})$ of non-intersecting lattice paths, where
$P_i$ runs from $(-i,-i)$ to $(i+\be,i+\be)$ and
does not pass over the diagonal $x=y$, $i=0,1,\dots,n-1$.
In \cite{\DeViAB}, Desainte--Catherine and Viennot showed  
that this counting problem is 
equivalent to the problem of counting tableaux with a bounded number 
of columns, all rows being of even length.  
Desainte--Catherine and Viennot solve the counting problem by
evaluating the determinant by means of 
the quotient-difference algorithm (see also \cite{\VienAB}). 
A weighted version of the tableaux counting problem of 
Desainte--Catherine and Viennot was solved  
by D\'es\-ar\-m\'enien \cite{\DesaAB, Th\'eor\`eme~1.2}.  
The determinant and its associated counting problem arise also
in the context of the determination of the multiplicity of Pfaffian
rings, see Theorem~2 and the accompanying remarks in \cite{\GhKrAA}.
\endremark

It is a simple observation that, given a determinant in which each
entry is the sum of two expressions,
$$\det_{0\le i,j\le n-1}(a_{i,j}+b_{i,j})$$
say, one can use linearity of the determinant in the rows to expand
this determinant into the sum
$$\sum _{S\subseteq \{0,1,\dots,n-1\}} ^{}
\det_{0\le i,j\le n-1}\(c^{(S)}_{i,j}\),
\tag\AI$$
where $c^{(S)}_{i,j}=a_{i,j}$ if $i\in S$ and $c^{(S)}_{i,j}=b_{i,j}$
otherwise. This sum consists of $2^n$ terms and, normally, 
will therefore not be very useful. Should it happen, however, that
the $(i+1)$-st row of the matrix $(a_{i,j})_{0\le i,j\le n-1}$ and the 
$i$-th row of the matrix $(b_{i,j})_{0\le i,j\le n-1}$ agree for all 
$i$, $i=0,1,\dots,n-2$, then
most terms in the sum (\AI) would vanish because, for a given set $S$,
there will always be two identical rows in the determinant
$\det_{0\le i,j\le n-1}\(c^{(S)}_{i,j}\)$ if there is an 
$s\in\{0,1,\dots,n-1\}$ with $s\notin S$ and $s+1\in S$. 
We summarise this last observation in the following lemma.

\proclaim{Lemma~\TC}
For all positive integers $n$ we have
$$\det_{0\le i,j\le n-1}(a_{i,j}+a_{i+1,j})=
\sum _{s=0} ^{n}
\det_{0\le i,j\le n-1}(a_{i+\chi(i\ge s),j}),
$$
where $\chi(\Cal S)=1$ if $\Cal S$ is
true and $\chi(\Cal S)=0$ otherwise.
\endproclaim

By combining Theorem~\TB\ and the above lemma, we obtain
the following corollary.

\proclaim{Corollary~\TD}
Let $n$ be a positive integer and $\al_0,\al_1,\dots,\al_{n}$
non-negative integers. Then
$$\multline 
\det_{0\le i,j\le n-1}\big(C_{\al_i+j}+C_{\al_{i+1}+j}\big)=
\prod _{0\le i<j\le n} ^{}(\al_j-\al_i)
\prod _{i=0} ^{n-1}\frac {(i+n)!} {(2i)!}
\prod _{i=0} ^{n}\frac {(2\al_i)!} {\al_i!\,(\al_i+n)!}\\
\times
\sum _{s=0} ^{n}\frac {\al_s!\,(\al_s+n)!} {(2\al_s)!\prod _{j=0}
^{s-1}(\al_s-\al_j)\prod _{j=s+1} ^{n}(\al_j-\al_s)}.
\endmultline$$
\endproclaim

\remark{Remarks}
(a) It can be checked that for $\al_i=i$, $i=0,1,\dots,n$,
Corollary~\TD\ reduces to
$$
\det_{0\le i,j\le n-1}\big(C_{i+j}+C_{i+j+1}\big)=
\sum _{s=0} ^{n}\binom {n+s}{n-s}=\sum _{s=0} ^{n}\binom {2n-s}s.
$$
In view of the well-known formula 
$$F_m=\sum _{s=0} ^{m}\binom {m-s}s$$
for Fibonacci numbers $F_m$, this immediately implies (\AE).

\smallskip
(b) Proofs of (\AE) by using non-intersecting lattice paths can be
found in \cite{\BeCQAA} and \cite{\CiglAS, Sec.~3}. 
\endremark

\subhead 4. Determinants of generalised Catalan numbers, I\endsubhead
The purpose of this section is to show how the determinant evaluation in
Lemma~\TA\ and Theorem~\TAA\ on non-intersecting lattice paths lead
to proofs of (\AF) and (\AG), and, in fact, of generalisations
thereof.

Let $\mu$ be a positive integer.
We denote by $P\big((0,0)\to (c,d)\mid x\ge \mu y\big)$ the number of
lattice paths starting at the origin and ending at $(c,d)$, which
never pass above the line $x=\mu y$. Then it is well-known (see
\cite{\MohaAE, Theorem~3}) that
$$
P\big((0,0)\to (c,d)\mid x\ge \mu y\big)=\frac {c-\mu d+1} {c+d+1}\binom
{c+d+1} {d}  .
\tag\AJ$$
The generalised Catalan numbers in (\AA) are a special case:
by (\AJ), the number $C_{n,k}$ is equal to the number of 
lattice paths starting at the origin and ending at
$\(n,\lfloor{\frac {n} {k-1}}\rfloor\)$, which never pass above the line $x=(k-1)y$.

\proclaim{Theorem~\TE}
Let $n$ and $k$ be positive integers and $\be,\al_0,\al_1,\dots,\al_{n-1}$
non-negative integers with $0\le \be\le k-1$. Then
$$
\det_{0\le i,j\le n-1}\big(C_{(k-1)\al_i+j+\be,k}\big)
=
\prod _{0\le i<j\le n-1} ^{}(\al_j-\al_i)
\prod _{i=0} ^{n-1}\frac {((k-1)i+\be+n)!\,(k\al_i+\be)!}
{(ki+\be)!\,\al_i!\,((k-1)\al_i+\be+n)!}.
$$
\endproclaim

\demo{Proof} 
By (\AJ) and Theorem~\TAA, we can
interpret the determinant in the theorem as the number of families
$(P_0,P_1,\dots,P_{n-1})$ of non-intersecting lattice paths, where
$P_i$ runs from $(-(k-1)\al_i,-\al_i)$ to $\(i+\be,\lfloor{\frac {i+\be}
{k-1}}\rfloor\)$ and
does not pass over the line $x=(k-1)y$, $i=0,1,\dots,n-1$. 

\midinsert
\vskip10pt
\vbox{
$$
\Gitter(10,6)(-7,-4)
\Koordinatenachsen(10,6)(-7,-4)
\Pfad(-6,-3),111111111111112222222\endPfad
\Pfad(-4,-2),1111111112112222\endPfad
\Pfad(-2,-1),111112112122\endPfad
\Pfad(0,0),1121121\endPfad
\PfadDicke{3pt}
\kern1pt
\Pfad(6,2),2\endPfad
\Pfad(7,2),2\endPfad
\Pfad(8,2),22\endPfad
\kern-1pt
\DickPunkt(-6,-3)
\DickPunkt(-4,-2)
\DickPunkt(-2,-1)
\DickPunkt(0,0)
\DickPunkt(5,2)
\DickPunkt(6,3)
\DickPunkt(7,3)
\DickPunkt(8,4)
\Kreis(6,2)
\Kreis(7,2)
\Kreis(8,2)
\Label\lo{A_0}(0,0)
\Label\lo{A_1}(-2,-1)
\Label\lo{A_2}(-4,-2)
\Label\lo{A_3}(-6,-3)
\Label\o{E_0}(5,2)
\Label\o{E_1}(6,3)
\Label\o{E_2}(7,3)
\Label\o{E_3}(8,4)
\Label\lo{P_0}(2,1)
\Label\r{P_1}(3,-1)
\Label\r{P_2}(5,-2)
\Label\ru{P_3}(8,-3)
\thinlines
\catcode`\@=11
\raise-4 \Einheit\hbox to0pt{\hskip-8 \Einheit
         \Line@(2,1){18}\hss}
\catcode`\@=12
\hskip1cm
$$
\centerline{\eightpoint A family of non-intersecting lattice paths}
\vskip7pt
\centerline{\eightpoint Figure \FA}
}
\endinsert

If we would directly apply Theorem~\TAA\ together with
the formula on the right-hand side of (\AJ) to the above problem of
counting non-intersecting lattice paths, then we would have to deal
with a determinant which is not easy to handle directly.
However, we may use combinatorics to simplify the determinant.
(The following simplifications could also be achieved by row and
column manipulations of the determinant. However, they become much
simpler and much more transparent in terms of non-intersecting
lattice paths.) The simplification is best explained with an example
at hand. Let us consider the case $k=3$, $n=4$, $\al_0=0$, $\al_1=1$, 
$\al_2=2$, $\al_3=3$, $\be=5$. Then the
starting and end points $A_0,\dots,A_3,E_0,\dots,E_3$ are as shown in 
Figure~\FA, which at the same time shows an example of a family of
non-intersecting lattice paths with these starting and end points.
(That portions of some paths are indicated by thick lines and some
points are circled should be ignored for the moment.)
Since the path family is non-intersecting, and since the paths
$P_1,P_2,\dots,P_{n-1}$ must stay (weakly) below the line $x=(k-1)y$ 
and must avoid $E_0$, they must pass to the right of $E_0$.
Hence, they must all end with vertical steps
above the height at which we find $E_0$, that is, above $y=\lfloor{\frac
{\be} {k-1}}\rfloor$. These vertical steps (in our example in Figure~\FA\
they are indicated by thick lines) can therefore be deleted without
changing the enumeration problem. Thus, we may equivalently count
families $P'_0,P'_1,\dots,P'_{n-1}$ of non-intersecting lattice
paths, where
$P'_i$ runs from\linebreak $A_i=(-(k-1)\al_i,-\al_i)$ to 
$E'_i=\(i+\be,\lfloor{\frac {\be} {k-1}}\rfloor\)$ and
does not pass over the line $x=(k-1)y$, $i=0,1,\dots,n-1$.
The modified family of paths which corresponds to our example in
Figure~\FA\ is shown in Figure~\FB.

\midinsert
\vskip10pt
\vbox{
$$
\Gitter(10,6)(-7,-4)
\Koordinatenachsen(10,6)(-7,-4)
\Pfad(-6,-3),1111111111111122222\endPfad
\Pfad(-4,-2),111111111211222\endPfad
\Pfad(-2,-1),11111211212\endPfad
\Pfad(0,0),1121121\endPfad
\DickPunkt(-6,-3)
\DickPunkt(-4,-2)
\DickPunkt(-2,-1)
\DickPunkt(0,0)
\DickPunkt(5,2)
\DickPunkt(6,2)
\DickPunkt(7,2)
\DickPunkt(8,2)
\Label\lo{A_0}(0,0)
\Label\lo{A_1}(-2,-1)
\Label\lo{A_2}(-4,-2)
\Label\lo{A_3}(-6,-3)
\Label\o{E'_0}(5,2)
\Label\o{E'_1}(6,2)
\Label\o{E'_2}(7,2)
\Label\o{E'_3}(8,2)
\Label\lo{P'_0}(2,1)
\Label\r{P'_1}(3,-1)
\Label\r{P'_2}(5,-2)
\Label\ru{P'_3}(8,-3)
\thinlines
\catcode`\@=11
\raise-4 \Einheit\hbox to0pt{\hskip-8 \Einheit
         \Line@(2,1){18}\hss}
\catcode`\@=12
\hskip1cm
$$
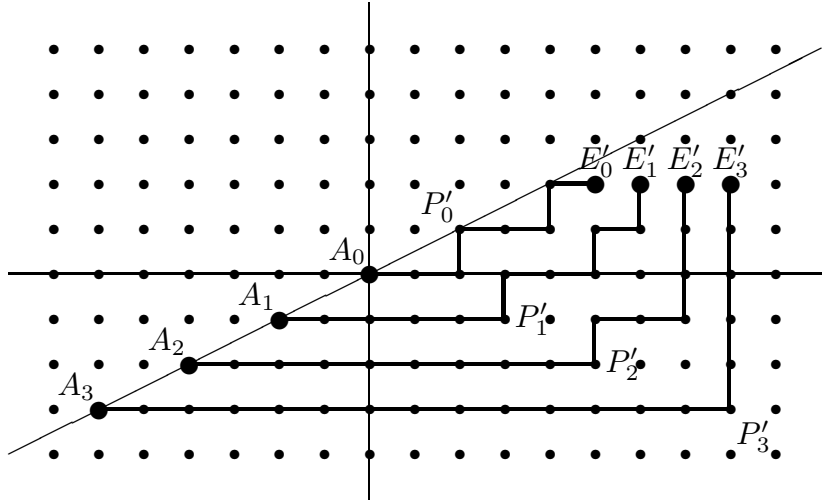
\centerline{\eightpoint The non-intersecting lattice paths of Figure~\FA\ 
without forced vertical steps}
\vskip7pt
\centerline{\eightpoint Figure \FB}
}
\endinsert

If we now use Theorem~\TAA\ and (\AJ),
then we obtain
$$
\det_{0\le i,j\le n-1}\big(C_{(k-1)\al_i+j+\be,k}\big)=
\det_{0\le i,j\le n-1}\(\frac {j+\be+1} {k\al_i+j+\be+1}
\binom {k\al_i+j+\be+1} {\al_i}\).
\tag\AJb
$$
By taking some factors out of the $i$-th row of the matrix of which
the determinant on the right-hand side is taken, $i=0,1,\dots,n-1$, we obtain
$$\align 
&\det_{0\le i,j\le n-1}\big(C_{(k-1)\al_i+j+\be,k}\big)=
\frac {(\be+n)!} {\be!}\prod _{i=0} ^{n-1}\frac {(k\al_i+\be)!}
{\al_i!\,((k-1)\al_i+\be+n)!}\\
&\kern1cm
\times
\det_{0\le i,j\le n-1}\big(
(k\al_i+\be+1)(k\al_i+\be+2)\cdots(k\al_i+\be+j)\\
&\kern2cm
\cdot
((k-1)\al_i+\be+j+2)((k-1)\al_i+\be+j+3)\cdots((k-1)\al_i+\be+n)
\big)\\
&\kern.5cm
=k^{\binom n2}(k-1)^{\binom n2}
\frac {(\be+n)!} {\be!}\prod _{i=0} ^{n-1}\frac {(k\al_i+\be)!}
{\al_i!\,((k-1)\al_i+\be+n)!}\\
&\kern3cm
\times
\det_{0\le i,j\le n-1}\bigg(
\(\al_i+\frac {\be+1}k\)\(\al_i+\frac {\be+2}k\)\cdots\(\al_i+\frac
{\be+j}k\)\\
&\kern4cm
\cdot
\(\al_i+\frac {\be+j+2} {k-1}\)
\(\al_i+\frac {\be+j+3} {k-1}\)\cdots\(\al_i+\frac {\be+n} {k-1}\)
\bigg).
\endalign$$
Now Lemma~\TA\ can be applied with $X_i=\al_i$, $A_j=(\be+j+1)/(k-1)$, and
$B_j=(\be+j)/k$. After some manipulation, one arrives at the
claimed result.\quad \quad \qed
\enddemo

\remark{Remark}
Again, from a conceptual point of view, 
the ``actual" theorem is the determinant evaluation for the
determinant on the right-hand side of (\AJb). 
In an equivalent form, this is
$$\multline 
\det_{0\le i,j\le n-1}\(
\binom {k\al_i+j+\be} {\al_i-1}\)
\\=\frac {\be!} {(\be+n)!}
\prod _{0\le i<j\le n-1} ^{}(\al_j-\al_i)
\prod _{i=0} ^{n-1}\frac {((k-1)i+\be+n)!\,(k\al_i+\be)!}
{(ki+\be)!\,(\al_i-1)!\,((k-1)\al_i+\be+n)!}.
\endmultline
\tag\AJa$$
By observing that 
$$\binom {k\al_i+j+\be} {\al_i-1}=
\binom {k\al_i+j+\be} {(k-1)\al_i+j+\be+1}=
(-1)^{(k-1)\al_i+j+\be-1}\binom {-\al_i} {(k-1)\al_i+j+\be-1},
$$
one sees that the determinant evaluation (\AJa) is equivalent to
\cite{\KratBN, Theorem~26, (3.13)}.
In particular, for the validity of (\AJa),
the restriction on $\be$ in the statement of Theorem~\TE\
is not necessary. Moreover, the identity (\AJa)
continues to hold for real or complex $\be$ and $\al_i$ if one
interprets the appearing factorials and binomials as suitable 
expressions involving gamma functions
(cf.\ \cite{\GrKPAA, \S5.5, (5.96), (5.100)}),
as long as there do not appear any singularities.
\endremark

By combining Theorem~\TE\ and Lemma~\TC, we obtain
the following corollary.

\proclaim{Corollary~\TF}
Let $n$ and $k$ be a positive integers and $\be,\al_0,\al_1,\dots,\al_{n}$
non-negative integers with $0\le \be\le k-1$. Then
$$\multline 
\det_{0\le i,j\le n-1}\big(C_{(k-1)\al_i+j+\be,k}+
   C_{(k-1)\al_{i+1}+j+\be,k}\big)\\
=
\prod _{0\le i<j\le n} ^{}(\al_j-\al_i)
\prod _{i=0} ^{n-1}\frac {((k-1)i+\be+n)!}
{(ki+\be)!}
\prod _{i=0} ^{n}\frac {(k\al_i+\be)!}
{\al_i!\,((k-1)\al_i+\be+n)!}\\
\times
\sum _{s=0} ^{n}\frac {\al_s!\,((k-1)\al_s+\be+n)!}
{(k\al_s+\be)!\prod _{j=0}
^{s-1}(\al_s-\al_j)\prod _{j=s+1} ^{n}(\al_j-\al_s)}.
\endmultline$$
\endproclaim

\remark{Remark} 
It can be checked that for $\al_i=i$, $i=0,1,\dots,n$,
Corollary~\TF\ reduces to (\AF).
\endremark

Finally, we address the determinant evaluation (\AG).
We shall actually consider the more general determinant
$$\det_{0\le i,j\le n-1}(C_{(k-1)i+j+\be,k}+C_{(k-1)i+j+\be+1,k}),
\tag\AK$$
where $\be\le k-1$.
Clearly, we can again apply Lemma~\TC\ (with the roles of rows and
columns interchanged) to obtain
$$\det_{0\le i,j\le n-1}(C_{(k-1)i+j+\be,k}+C_{(k-1)i+j+\be+1,k})=
\sum _{s=0} ^{n}
\det_{0\le i,j\le n-1}(C_{(k-1)i+j+\chi(j\ge s)+\be,k}).
\tag\AL$$
By Theorem~\TAA, for any fixed $s$,
the determinant on the right-hand side has a
combinatorial interpretation in terms of non-intersecting lattice
paths in a manner completely analogous to the one in the proof of
Theorem~\TE: it counts families 
$(P_0,P_1,\dots,P_{n-1})$ of non-intersecting lattice paths, where
$P_i$ runs from $(-(k-1)i,-i)$ to 
$\big(i+\chi(i\ge s)+\be,\mathbreak
\lfloor{\frac {i+\chi(i\ge s)+\be} {k-1}}\rfloor\big)$ and
does not pass over the line $x=(k-1)y$, $i=0,1,\dots,n-1$.
Figure~\FC\ shows an example with $n=4$, $k=3$, $\be=1$ and $s=2$.

\midinsert
\vskip10pt
\vbox{
$$
\Gitter(7,4)(-7,-4)
\Koordinatenachsen(7,4)(-7,-4)
\Pfad(-6,-3),1111111111212222\endPfad
\Pfad(-4,-2),111111122122\endPfad
\Pfad(-2,-1),111122\endPfad
\Pfad(0,0),1\endPfad
\PfadDicke{3pt}
\Pfad(-6,-3),111111111\endPfad
\Pfad(-4,-2),1111111\endPfad
\Pfad(-2,-1),111122\endPfad
\kern1pt
\Pfad(0,0),1\endPfad
\Pfad(4,1),2\endPfad
\Pfad(5,1),2\endPfad
\kern-1pt
\DickPunkt(-6,-3)
\DickPunkt(-4,-2)
\DickPunkt(-2,-1)
\DickPunkt(0,0)
\DickPunkt(1,0)
\DickPunkt(2,1)
\DickPunkt(4,2)
\DickPunkt(5,2)
\Kreis(3,-3)
\Kreis(3,-2)
\Kreis(4,1)
\Kreis(5,1)
\Label\lo{A_0}(0,0)
\Label\lo{A_1}(-2,-1)
\Label\lo{A_2}(-4,-2)
\Label\lo{A_3}(-6,-3)
\Label\o{E_0}(1,0)
\Label\o{E_1}(2,1)
\Label\o{E_2}(4,2)
\Label\o{E_3}(5,2)
\Label\ru{\raise4pt\hbox{$P_0$}}(0,0)
\Label\r{\hbox{$\kern-5pt P_1$}}(2,0)
\Label\r{P_2}(3,-1)
\Label\ru{P_3}(4,-3)
\thinlines
\catcode`\@=11
\raise-4 \Einheit\hbox to0pt{\hskip-8 \Einheit
         \Line@(2,1){15}\hss}
\catcode`\@=12
\hskip-1cm
$$
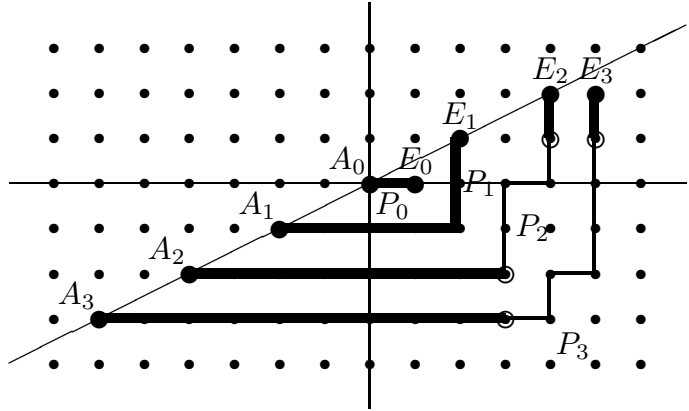
\centerline{\eightpoint Another family of non-intersecting lattice paths}
\vskip7pt
\centerline{\eightpoint Figure \FC}
}
\endinsert

Again, if the paths want to stay below the line $x=(k-1)y$ (being allowed
to touch it) and to be non-intersecting, then there are forced path
portions. These come in three different flavours. First, for $0\le
i<s$, the path $P_i$ must start with $ki+\be$ horizontal steps,
followed by as many vertical steps as necessary to reach the end point
$E_i$. (See the paths $P_0$ and $P_1$ in Figure~\FC.) Second, 
for $s\le i\le n-1$, the path $P_i$ must start with $(k-1)i+s+\be$
horizontal steps. (See the thick portions of $P_2$ and $P_3$ below
the $x$-axis in Figure~\FC.) Third, the paths $P_s,P_{s+1},\dots,P_{n-1}$ 
must all end with vertical steps
above the height at which we find $E_s$, that is, above $y=\lfloor{\frac
{s+\be+1} {k-1}}\rfloor$. Moreover, if $E_s$ should lie on the line $x=(k-1)y$
(that is, if $k-1$ divides $s+\be+1$),
then all of $P_s,P_{s+1},\dots,P_{n-1}$ must have an additional
vertical step from height $\frac {s+\be+1} {k-1}-1$ to $\frac
{s+\be+1} {k-1}$.
(See the thick vertical parts of $P_2$ and $P_3$ in
Figure~\FC. In the example, we have indeed that $E_s=E_2$ lies on
$x=(k-1)y=2y$.) The last observation can also be succinctly rephrased by 
saying that the paths $P_s,P_{s+1},\dots,P_{n-1}$ 
must all end with vertical steps above $y=\lfloor{\frac {s+\be}
{k-1}}\rfloor$.

\midinsert
\vskip10pt
\vbox{
$$
\Gitter(7,4)(-7,-4)
\Koordinatenachsen(7,4)(-7,-4)
\Pfad(3,-3),121222\endPfad
\Pfad(3,-2),2212\endPfad
\DickPunkt(3,-3)
\DickPunkt(3,-2)
\DickPunkt(4,1)
\DickPunkt(5,1)
\Label\l{A'_2}(3,-2)
\Label\l{A'_3}(3,-3)
\Label\o{E'_2}(4,1)
\Label\o{E'_3}(5,1)
\Label\ru{P'_2}(3,0)
\Label\ru{P'_3}(4,-3)
\thinlines
\catcode`\@=11
\raise-4 \Einheit\hbox to0pt{\hskip-8 \Einheit
         \Line@(2,1){15}\hss}
\catcode`\@=12
\hskip-1cm
$$
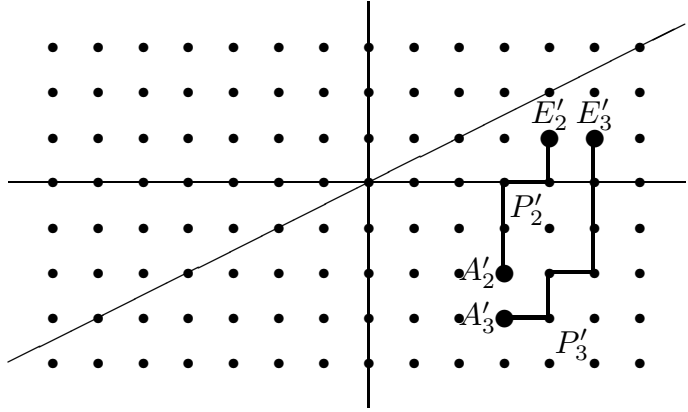
\centerline{\eightpoint The non-intersecting lattice paths of
Figure~\FC\ after removing forced path portions}
\vskip7pt
\centerline{\eightpoint Figure \FD}
}
\endinsert

We may delete the forced portions without changing the enumeration
problem. (See Figure~\FD\ for the resulting path family after the
forced portions have been removed in Figure~\FC.)
The boundary $x=(k-1)y$ may also be removed without changing
the enumeration problem. (See Figure~\FE\ for the result in case
of our running example from Figures~\FC\ and \FD. The dotted path
should be ignored at the moment.)
In that manner, we see that the determinant
on the right-hand side of (\AL) is equal to the number of families 
$(P'_s,P'_{s+1},\dots,P'_{n-1})$ of non-intersecting lattice paths, where
$P'_i$ runs from $(s+\be,-i)$ to 
$\(i+\be+1,\lfloor{\frac {s+\be} {k-1}}\rfloor\)$, $i=s,s+1,\dots,n-1$.

\midinsert
\vskip10pt
\vbox{
$$
\Gitter(7,4)(0,-4)
\Koordinatenachsen(7,4)(0,-4)
\Pfad(3,-3),121222\endPfad
\Pfad(3,-2),2212\endPfad
\SPfad(3,1),64646\endSPfad
\DickPunkt(3,-3)
\DickPunkt(3,-2)
\DickPunkt(4,1)
\DickPunkt(5,1)
\Kreis(3,1)
\Kreis(5,-4)
\Label\l{A'_2}(3,-2)
\Label\l{A'_3}(3,-3)
\Label\o{E'_2}(4,1)
\Label\o{E'_3}(5,1)
\Label\lo{S}(3,1)
\Label\ru{T}(5,-4)
\hskip4cm
$$
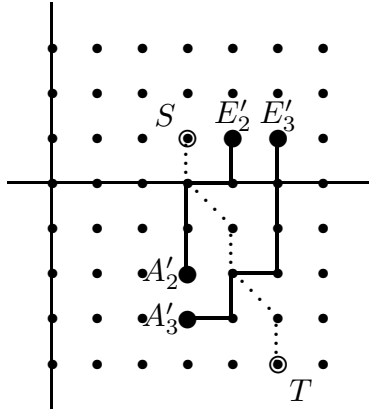
\centerline{\eightpoint The dual path for the non-intersecting
lattice paths of Figure~\FD}
\vskip7pt
\centerline{\eightpoint Figure \FE}
}
\endinsert

A more general enumeration problem can actually be solved in closed form.
We formulate the corresponding result in the proposition below.

\proclaim{Proposition~\TG}
Let $a,b,c$ and $\al_0,\al_1,\dots,\al_{n-1}$ be integers with
$a\le \al_0\le\al_1\le\dots\le\al_{n-1}$ and $b\le c$.
Then the number of families $(P_0,P_1,\dots,P_{n-1})$ of non-intersecting
lattice paths, the path $P_i$ running from $(a,b-i)$ to $(\al_i,c)$,
$i=0,1,\dots,n-1$, is given by
$$
\prod _{0\le i<j\le n-1} ^{}(\al_j-\al_i)
\prod _{i=0} ^{n-1}\frac {(\al_i+c-a-b)!} {(\al_i-a)!\,(c-b+i)!}.
$$
\endproclaim

\demo{Proof}
By Theorem~\TAA, the number in
question is equal to the determinant
$$\det_{0\le i,j\le n-1}\(\binom {\al_j-a+c-b+i}{c-b+i}\).$$ 
By taking some factors out of the $j$-th column of the matrix of which
we want to compute the determinant, $j=0,1,\dots,n-1$, 
we obtain the equivalent expression
$$\multline 
\prod _{j=0} ^{n-1}\frac {(\al_j+c-a-b)!} {(\al_j-a)!\,(c-b+j)!}\\
\times
\det_{0\le i,j\le n-1}\big((\al_j+c-a-b+1)(\al_j+c-a-b+2)\cdots
(\al_j+c-a-b+i)\big).
\endmultline$$
This determinant is of the form $\det_{0\le i,j\le
n-1}\big(p_i(\al_j)\big)$, where $p_i(x)$ is a polynomial in $x$ of
degree $i$ with leading coefficient $1$, $i=0,1,\dots,n-1$. 
Such a determinant can be easily
reduced to the Vandermonde determinant $\det_{0\le i,j\le
n-1}\big(\al_j^i\big)$ by using elementary column operations
(see also \cite{\KratBN, Proposition~1}). The evaluation of the
Vandermonde determinant being well-known, the assertion of our
proposition follows immediately.\quad \quad \qed
\enddemo

We now apply Proposition~\TG\ to obtain the announced generalisation
of (\AG).

\proclaim{Theorem~\TH}
Let $n$ and $k$ be positive integers and $\be$ a
non-negative integer with $0\le \be\le k-1$. Then
$$
\det_{0\le i,j\le
n-1}\big(C_{(k-1)i+j+\be,k}+C_{(k-1)i+j+\be+1,k}\big)
=
\sum _{s=0} ^{n}\binom {\lfloor \frac {s+\be} {k-1}\rfloor+n} {n-s}.
\tag\ALa$$
\endproclaim

\demo{Proof}
Coming back to the enumeration problem which is illustrated in
Figure~\FE, we apply Proposition~\TG\ with $n$ replaced by $n-s$,
$a=s+\be$, $b=-s$, $c=\lfloor\frac {s+\be} {k-1}\rfloor$, and
$\al_i=s+i+\be+1$, $i=0,1,\dots,n-s-1$.
The claimed result follows upon little simplification.\quad \quad \qed
\enddemo

\remark{Remarks}
(a) Clearly, Eq.~(\AG) is the special case $\be=0$ of Theorem~\TH.

\smallskip
(b) It should be clear that, from a conceptual point of view, 
Proposition~\TG\ represents the actual key result which is behind 
Theorem~\TH. On the other hand, if we are just interested in proving
Theorem~\TH, then we do not need the determinant evaluation
in the proof of Proposition~\TG. Namely, for Theorem~\TH, we only
need the special case of Proposition~\TG\ in which the end points of
the paths $P_i$ are successive lattice points on a horizontal line,
and the first end point is by one unit to the right of the vertical
line on which we find the starting points (see Figure~\FE).
The reader is referred to the proof of Theorem~\TH\ for the exact
choices of the starting and end points.

There is an elegant argument using the concept of {\it dual paths},
an idea which is again due to Gessel and Viennot \cite{\GeViAA,
Sec.~4}, to see that the number of families of non-intersecting
lattice paths with starting and end points as described above is a
binomial coefficient. In order to explain the idea, we let $S$ be the
lattice point one unit to the left of the first end point,
that is, $S=(s+\be,\lfloor{\frac {s+\be} {k-1}}\rfloor)$, and we let
$T$ be the lattice point exactly below the last end point, the height 
of which is by one unit lower than the height of the last starting point, 
that is, $T=(n+\be,-n)$. See Figure~\FE.
We now connect $S$ with $T$ by moving vertically 
downwards, unless we hit one of the existing paths. If the
latter happens, then we continue by a diagonal step $(1,-1)$, etc.
The ``dual" path which results for our example in Figure~\FE\ is
indicated by dotted line segments.

It is easy to see that families of non-intersecting lattice paths
connecting the starting and end points as above are in bijection with
all paths from $S$ to $T$ consisting of vertical down steps $(0,-1)$
and diagonal down steps $(1,-1)$. However, each such path from $S$ to
$T$ has exactly $s+\lfloor{\frac {s+\be} {k-1}}\rfloor$ vertical steps and 
$n-s$ diagonal steps, the order of which can be chosen freely. Thus,
there are $\binom {\lfloor \frac {s+\be} {k-1}\rfloor+n} {n-s}$
such paths and, therefore, as many families of non-intersecting
lattice paths with the starting and end points as above, in
accordance with (\ALa).

For a slightly different lattice path proof of (\AA) see \cite{\BeCQAB}. 
\endremark

\subhead 5. Determinants of generalised Catalan numbers, II\endsubhead
Let again $k$ be a fixed positive integer.
Rewriting the Catalan number $C_n$ in the form
$$C_n=\frac {1} {2n+1}\binom {2n+1}{n+1}=\frac {1} {2n+1}\binom
{2n+1}{n},\tag\AM$$
another way to construct ``generalised" Catalan numbers is by
considering numbers of the form
$$\frac {l} {kn+l}\binom {kn+l}{n}.$$
Tamm \cite{\TammAD} has studied Hankel determinants involving such
numbers to considerable extent. It seems that, beyond $k=3$, one cannot
expect any closed form results
(see also \cite{\EgRRAB}). On the other hand, for $k=3$ there
exist several beautiful results, with many variations, due to 
E\u gecio\u glu, Redmond and Ryavec \cite{\EgRRAA} (for some of them) 
and Gessel and Xin \cite{\GeXiAA} (for most of them). 
As the proof methods show, these results
lie on a much deeper level than their counterparts in Theorems~\TB\ or \TE.
In the theorem below, we present the determinant evaluations
involving generalised Catalan numbers of the form (\AM). 
For a complete list see \cite{\KratBZ, Theorem~31}.

\proclaim{Theorem~\TI} 
For any positive integer $n$, there hold 
$$ 
\det_{0\le i,j\le n-1}\(\frac {1} {3i+3j+1}\binom {3i+3j+1}{i+j}\) 
=
\prod _{i=0} ^{n-1}\frac {(\frac {2} {3})_i\,(\frac {1} {6})_i\,
(\frac {4} {3})_i\,(\frac {5} {6})_i} 
{(\frac {1} {2})_{2i}\,(\frac {3} {2})_{2i}}\(\frac {27} {4}\)^{2i},
\tag\AN
$$
$$
\det_{0\le i,j\le n-1}\(\frac {1} {3i+3j+4}\binom {3i+3j+4}{i+j+1}\) 
=
\prod _{i=0} ^{n-1}\frac {(\frac {4} {3})_i\,(\frac {5} {6})_i\,
(\frac {5} {3})_i\,(\frac {7} {6})_i} 
{(\frac {3} {2})_{2i}\,(\frac {5} {2})_{2i}}\(\frac {27} {4}\)^{2i},
\tag\AO
$$
$$
\det_{0\le i,j\le n-1}\(\frac {1} {3i+3j+2}\binom {3i+3j+2}{i+j+1}\) 
=
\prod _{i=0} ^{n-1}\frac {(\frac {4} {3})_i\,(\frac {5} {6})_i\,
(\frac {5} {3})_i\,(\frac {7} {6})_i} 
{(\frac {3} {2})_{2i}\,(\frac {5} {2})_{2i}}\(\frac {27} {4}\)^{2i},
\tag\AP
$$
$$
\det_{0\le i,j\le n-1}\(\frac {1} {3i+3j+5}\binom {3i+3j+5}{i+j+2}\) 
=
\prod _{i=0} ^{n}\frac {(\frac {2} {3})_i\,(\frac {1} {6})_i\,
(\frac {4} {3})_i\,(\frac {5} {6})_i} 
{(\frac {1} {2})_{2i}\,(\frac {3} {2})_{2i}}\(\frac {27} {4}\)^{2i},
\tag\AQ
$$
$$
\det_{0\le i,j\le n-1}\(\frac {2} {3i+3j+1}\binom {3i+3j+1}{i+j+1}\) 
=
\prod _{i=0} ^{n}\frac {(\frac {2} {3})_i\,(\frac {1} {6})_i\,
(\frac {4} {3})_i\,(\frac {5} {6})_i} 
{(\frac {1} {2})_{2i}\,(\frac {3} {2})_{2i}}\(\frac {27} {4}\)^{2i},
\tag\AR
$$
$$
\det_{0\le i,j\le n-1}\(\frac {2} {3i+3j+4}\binom {3i+3j+4}{i+j+2}\) 
=
\prod _{i=0} ^{n}\frac {(\frac {4} {3})_i\,(\frac {5} {6})_i\,
(\frac {5} {3})_i\,(\frac {7} {6})_i} 
{(\frac {3} {2})_{2i}\,(\frac {5} {2})_{2i}}\(\frac {27} {4}\)^{2i},
\tag\AS
$$
where $(\al)_i:=\al(\al+1)\cdots(\al+i-1)$ is the usual Pochhammer symbol.
Let $a_0=-2$ and $a_m=\frac {1} {3m+1}\binom {3m+1}m$ for $m\ge1$. Then
$$ 
\det_{0\le i,j\le n-1}(a_{i+j})=
\prod _{i=0} ^{n-1}(-2)\frac {(\frac {1} {3})_i\,(-\frac {1} {6})_i\,
(\frac {5} {3})_i\,(\frac {7} {6})_i} 
{(\frac {1} {2})_{2i}\,(\frac {3} {2})_{2i}}\(\frac {27} {4}\)^{2i}.
\tag\AT
$$
Let $b_0=10$ and $b_m=\frac {2} {3m+2}\binom {3m+2}m$ for $m\ge1$. Then
$$ 
\det_{0\le i,j\le n-1}(b_{i+j})=
\prod _{i=0} ^{n-1}10\frac {(\frac {2} {3})_i\,(\frac {1} {6})_i\,
(\frac {7} {3})_i\,(\frac {11} {6})_i} 
{(\frac {3} {2})_{2i}\,(\frac {5} {2})_{2i}}\(\frac {27} {4}\)^{2i}.
\tag\AU
$$
Furthermore,
$$ 
\det_{0\le i,j\le n-1}\(\frac {2} {3i+3j+5}\binom {3i+3j+5}{i+j+1}\) 
=
\prod _{i=0} ^{n}\frac {(\frac {4} {3})_i\,(\frac {5} {6})_i\,
(\frac {5} {3})_i\,(\frac {7} {6})_i} 
{(\frac {3} {2})_{2i}\,(\frac {5} {2})_{2i}}\(\frac {27} {4}\)^{2i}.
\tag\AV
$$
Let $c_0=\frac {7} {2}$ and 
$c_m=\frac {2} {3m+1}\binom {3m+1}{m+1}$ for $m\ge1$. Then
$$ 
\det_{0\le i,j\le n-1}(c_{i+j})=
\prod _{i=0} ^{n}\frac {(\frac {4} {3})_i\,(\frac {5} {6})_i\,
(\frac {5} {3})_i\,(\frac {7} {6})_i} 
{(\frac {3} {2})_{2i}\,(\frac {5} {2})_{2i}}\(\frac {27} {4}\)^{2i}.
\tag\AW
$$
\endproclaim

\enddocument